\documentclass{amsart}
\usepackage[latin1]{inputenc}
\usepackage[english]{babel}
\usepackage{amsmath,amssymb}
\usepackage{float}
\IfFileExists{mathpazo.sty}{
\usepackage{mathpazo}
}{}
\usepackage{mathrsfs}
\usepackage[T1]{fontenc}
\usepackage{amscd}
\usepackage{boxedminipage}
\usepackage[line,cmtip,all,matrix,dvips]{xy}
\usepackage{graphicx}

\renewcommand{\epsilon}{\varepsilon}
\renewcommand{\rho}{\varrho}
\renewcommand{\phi}{\varphi}

\newtheorem{proposition}{Proposition}[section]
\newtheorem{lemma}{Lemma}[section]
\newtheorem{theorem}{Theorem}[section]
\newtheorem{corollary}{Corollary}[section]
\newtheorem{definition}{Definition}[section]

\newtheorem{remark}{Remark}[section]

\newcommand{\eproof}{\begin{flushright} $\square$ \end{flushright}}

\newcommand{\stable}{stable\mbox{ }}

\newcommand{\fm}{saturated\mbox{ }}

\newcommand{\e}[1]{\mathbf #1}

\newcommand{\A}{\mathcal A}

\newcommand{\tF}{{\tilde F}}

\newcommand{\Bosonnormalord}{\stackrel{\textstyle \circ}{ \circ}}

\newcommand{\spane}{\mathop{\fam0 Span}\nolimits}

\newcommand{\id}{\mathop{\fam0 Id}\nolimits}

\newcommand{\Id}{\mathop{\fam0 Id}\nolimits}

\newcommand{\R}{\mathop{\fam0 {\mathbb R}}\nolimits}

\newcommand{\C}{\mathop{\fam0 {\mathbb C}}\nolimits}

\newcommand{\Z}{{\mathbb Z}}

\newcommand{\ra}{\mathop{\fam0 \rightarrow}\nolimits}

\newcommand{\M}{\mathop{\fam0 {\mathcal M}}\nolimits}

\newcommand{\Kappa}{\Psi}

\newcommand{\tc}{\mathop{\fam0 {\tilde c}}\nolimits}
\newcommand{\tC}{\mathop{\fam0 {\tilde C}}\nolimits}
\newcommand{\tPi}{{\tilde \Pi}}

\def\l{\lambda}\def\m{\mu}
\def\L{\Lambda}

\newcommand{\Si}{\Sigma}
\newcommand{\Sib}{{\mathbf \Sigma}}
\newcommand{\Thb}{{\mathbf \Theta}}

\def\:{:}

\def\Bosonnormalordconstruction#1{\vcenter{\hbox{\ooalign{%
\raise.8ex\hbox{$#1\circ$}\crcr\lower.8ex\hbox{$#1\circ$}}}}}

\def\lg{\lambda_\gamma}
\newcommand{\ga}{\gamma}

\newcommand{\tP}{\tilde P}\newcommand{\tV}{\tilde V}
\newcommand{\tSib}{\tilde \Sib}\newcommand{\tSi}{\tilde \Si}

\def\Bosonnormalord{\,\lower.8ex \hbox{$\circ$} \llap{\raise.8ex\hbox{$\circ$}} \,}
\def\normalord{\,\lower.8ex \hbox{$\cdot$} \llap{\raise.8ex\hbox{$\cdot$}} \,}

\def\ra{\rightarrow}

\def\tL{{\tilde L}}

\def\bC{{\mathbf C}}

\def\bR{{\mathbf R}}
\def\bZ{{\mathbf Z}}

\def\cD{{\mathcal D}}
\def\cF{{\mathcal F}}
\def\cS{{\mathcal S}}

\def\tnu{{\tilde \nu}}

\def\Id{\mathop{\rm Id}\nolimits}

\begin{document}

\title{Modular functors are determined by their genus zero data.}

\author{J{\o}rgen Ellegaard Andersen}
\address{Department of Mathematics\\
        University of Aarhus\\
        DK-8000, Denmark}
\email{andersen{\@@}imf.au.dk}

\author{Kenji Ueno}
\address{Department of Mathematics\\
        Faculty of Science, Kyoto University\\
        Kyoto, 606-01 Japan}
\email{ueno{\@@}kusm.kyoto-u.ac.jp}

\begin{abstract}  We prove in this paper that the genus zero data of a modular
functor determines the modular functor. We do this by establishing
that the $S$-matrix in genus one with one point labeled arbitrarily
can be expressed in terms of the genus zero information and we give
an explicit formula. We do not assume the modular functor in
question has duality or is unitary, in order to establish this.
\end{abstract}
\maketitle

\tableofcontents

\section{Introduction}

A modular functor is a functor $Z$ from the category of smooth
surfaces with certain extra structure, namely the category of
labeled marked surfaces, to the category of finite dimensional
vector spaces over the complex numbers. In section \ref{AxiomsMF} we
shall give the precise axioms for a modular functor following K.
Walker's topological reformulation \cite{Walker} of G. Segal's
axioms for a conformal field theory \cite{Se}.

The objects of the category of labeled marked surfaces are pairs
consisting of a marked surface and a labeling. A marked surface here
refers to a quadruple of structures: A smooth closed oriented
surface, a finite set of "marked" points, a tangent direction at
each of the marked points together with a Lagrangian subspace of the
first homology of the surface. Labeled means that each marked point
of the surface is labeled by an element from a certain finite
label-set $\Lambda$ which is specific to the modular functor $Z$.
This label set is further required to have an involution $\mbox{
}^\dagger$ and a preferred $0\in\Lambda$, such that $0^\dagger = 0$.

By the factorization property of a modular functor, we can
express the vector space associated to any label marked surface
as a direct sum of tensor products of vector spaces associated to
2-spheres with three marked points labeled appropriately.
This is done by choosing a pair of pants decomposition of the
surface.

Let $S^2 = \C{} \cap \{\infty\}$ and $v_t$
be the direction along the positive real axis at $t\in \bR \cup
\infty \subset S^2$.
We let $\Upsilon = (S^2;0,1,\infty; v_0,v_1,v_\infty)$. For
$\lambda,\mu,\nu \in \Lambda$, we define
\[Z_{\lambda,\mu,\nu} = Z(\Upsilon, \lambda,\mu,\nu).\]
\begin{center}
\begin{figure}[H]
\includegraphics[scale=.5]{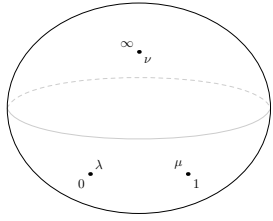}
\caption{A sphere with three labeled marked points.}
\label{figs-0}
\end{figure}
\end{center}

\begin{center}
\begin{figure}[H]
\includegraphics[scale=.5]{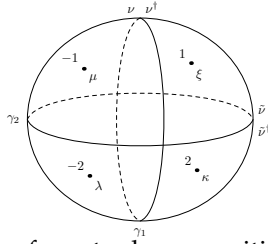}
\caption{Two pair of pants decompositions of a sphere with four labeled marked points.}
\label{figs-31}
\end{figure}
\end{center}

Any two pair of pants decompositions of a given surface
can be obtained, one from the other by a sequence of flips.
There
are two kinds of flips needed.
The first one is the change from
the pair of pants decomposition of the 2-sphere with four marked
points given by $\gamma_1$ to the one given by $\gamma_2$ in
figure \ref{figs-31}.

As will be explained in section \ref{BD}, this gives us an isomorphism
\[F\left[ \begin{array}{cc}
  \mu & \xi \\
  \lambda & \kappa
\end{array} \right] : \bigoplus_{\nu \in \Lambda} Z_{\nu, \mu, \lambda  }
\otimes Z_{\nu^\dagger, \kappa, \xi } \ra
\bigoplus_{\tnu \in \Lambda}Z_{\tnu, \lambda, \kappa }
\otimes Z_{\tnu^\dagger, \xi, \mu }.\]

\begin{center}
\begin{figure}[H]
\includegraphics[scale=.5]{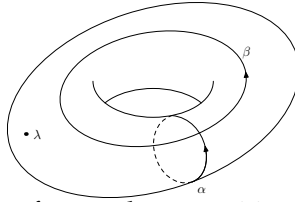}
\caption{Two pair of pants decompositions of a torus with one labeled marked point.}
\label{figs-3}
\end{figure}
\end{center}

 The second kind of flip is the change from the
pair of pants decomposition of the torus with one marked point
given by $\alpha$ to the one given by $\beta$ in figure
\ref{figs-3}.

For each label $\lambda\in \Lambda$, we get an isomorphism
$$S(\lambda) : \bigoplus_{\mu} Z_{\lambda, \mu,\mu^\dagger} \ra
\bigoplus_{\nu} Z_{\lambda, \nu,\nu^\dagger}.$$ For $\lambda = 0$,
we simply write $S = S(0)$. From the axioms of a modular functor
it follows that $\dim Z_{0,\mu,\mu^\dagger} = 1$ for all
$\mu\in\Lambda$. Further, we will see in section \ref{BD}, that
the axioms of a modular functor determines a unique non-zero
vector in $Z_{0,\mu,\mu^\dagger}$, and so we get a matrix $S =
(S_{\mu,\nu})_{\mu,\nu \in \Lambda}$.

A flip relating two pair of pants decompositions gives
an isomorphism between the corresponding two direct sums of
tensor products of the vector spaces $Z_{\lambda,\mu,\nu}$ for appropriately varying
$\lambda,\mu,\nu \in \Lambda$.

For the first flip, the corresponding isomorphism is clearly
determined by the functor applied to genus zero surfaces with less
than or equal to four marked points. On the other hand the second
flip involves a genus one surface with one marked point. However, we
will prove in this paper that $S(\lambda)$ for all $\lambda \in
\Lambda$ is determined by the restriction of the modular functor to
genus zero surfaces. We call this restriction to genus zero surfaces
of a modular functor the genus zero data of the modular functor.

\begin{theorem}\label{Main}
The genus zero data of any modular functor determines the modular
functor.
\end{theorem}

The key ingredient in the proof of this theorem is what we call the
curve operators.  Given an oriented
simple closed curve $\gamma$ on a labeled marked surface
and a label $\lambda_\gamma\in \Lambda$, we construct an
endomorphism of the vector space the
modular functor associates to the labeled marked surface.
  Loosely speaking they are obtained
 by creating two points labeled by $\lambda_\gamma$ and $\lambda^\dagger_\gamma$ near
 each other along $\gamma$. Move one of them around $\gamma$ and then
 annihilate them again. For the precise definition see sections \ref{CO}.
 The next step is to express
the automorphism induced by a Dehn-twist in a simple closed curve
as a linear combination of curve operators for the same
curve. In fact this linear combination is seen to be universal by factoring
along the boundary of a tubular neighbourhood of the curve. Moreover, it is then clear the linear combination
is also determined by genus zero data. Now, since curve operators are determined completely
by genus zero data, as we argue in section \ref{From1to0}, we see
that the Dehn-twist in any simple closed curve is determined by
genus zero data. Using the standard presentation of the mapping
class group of surface of genus one with one marked point, we
conclude that the matrices $S(\lambda)$ are also determined by the
genus zero data.

In fact, we will establish the following explicit formula for
$S(\lambda)$.

Pick a basis $\zeta_j(\lambda,\mu,\nu)$, $j= 1,\ldots \dim
Z_{\lambda,\mu,\nu}$ for $Z_{\lambda,\mu,\nu}$. For $\nu=
\mu^\dagger$ and $\lambda = 0$, we will assume that
$\zeta_1(0,\mu,\mu^\dagger)$ is this preferred vector in $Z_{0,\mu,\mu^\dagger}$, as
discussed
above. We then define
\begin{eqnarray*}
\lefteqn{
F_{\nu,\tnu}\left[ \begin{array}{cc}
  \mu & \xi \\
  \lambda & \kappa
\end{array} \right](\zeta_i(\nu, \mu,\lambda) \otimes
\zeta_j(\nu^\dagger,\kappa, \xi))}\\
& = & \sum_{k,l} F_{\nu,\tnu}\left[ \begin{array}{cc}
  \mu & \xi \\
  \lambda & \kappa
\end{array} \right]_{ij}^{kl}
\zeta_k(\tnu, \lambda,\kappa) \otimes \zeta_l(\tnu^\dagger,\xi,\mu)
\end{eqnarray*}

and

\[S(\lambda) (\zeta_i(\lambda,\mu,\mu^\dagger)) =
\sum_{\nu,j} S(\lambda)_{\mu,i}^{\nu,j}\zeta_j(\lambda,\nu,\nu^\dagger). \]

Also we need the following three self diffeomorphisms of $\Upsilon$.
The diffeomorphism $R : \Upsilon \ra \Upsilon$ is  $R(z) =
(z-1)/z$. It induces a linear isomorphism
\[R : Z_{\lambda,\mu,\nu} \ra Z_{\mu,\nu,\lambda}\]
with the matrix presentation
\[R(\zeta_i(\lambda,\mu,\nu)) = \sum_j R_{ij}\zeta_j(\mu,\nu,\lambda).\]

The diffeomorphism $B : \Upsilon \ra \Upsilon$ is the composition
of $z \mapsto z/(z-1)$ with a negative half-twist at $0$ and
$\infty$ and a positive half-twist at $1$. Again it induces a
linear isomorphism
\[B : Z_{\lambda,\mu,\nu} \ra Z_{\lambda,\nu,\mu}\]
with the matrix presentation
\[B(\zeta_i(\lambda,\mu,\nu)) = \sum_j B_{ij}\zeta_j(\lambda,\nu,\mu).\]

Finally, the Dehn twist in a circle centered in $0$ and
of radius bigger than $1$ induces a further endomorphism
of $\Upsilon$. In particular it
induces an isomorphism of $Z_{0,\mu,\mu^\dagger}$ to it self. But
since this space is one dimensional, this isomorphism is simply
given by multiplication by a non-zero complex number $d_\mu$.

We introduce the so call twisted $F$-isomorphism
$$
\tF_{\nu,\tnu}\left[ \begin{array}{cc}
  \mu & \xi \\
  \lambda & \kappa
\end{array} \right] : \bigoplus_{\nu \in \Lambda}
Z_{\nu,\mu,\lambda} \otimes Z_{\xi,\nu^\dagger,\kappa} \ra
\bigoplus_{\tnu \in \Lambda} Z_{\kappa,\lambda,\tnu}
\otimes Z_{\tnu^\dagger,\xi,\mu}$$
defined by
$$
\tF_{\nu,\tnu}\left[ \begin{array}{cc}
 \mu & \xi \\
  \lambda & \kappa
\end{array} \right] = (B R^2 \otimes \Id) F_{\nu,\tnu}\left[ \begin{array}{cc}
  \mu & \xi \\
  \lambda & \kappa
\end{array} \right] (\Id \otimes B R).$$
In the matrix presentation we get
$$
\tF_{\nu,\tnu}\left[ \begin{array}{cc}
   \mu & \xi \\
  \lambda & \kappa
\end{array} \right]_{ki}^{jm} = \sum_{p,r,s,w}R_{ip}B_{pr} R^2_{sw}B_{wj}
 F_{\nu,\tnu}\left[ \begin{array}{cc}
 \mu & \xi \\
  \lambda & \kappa
\end{array} \right]_{kr}^{sm}.$$

\begin{theorem}\label{Main'}
For any modular functor we have the following formula
\[S(\lambda)_{\mu,i}^{\nu,j} = \sum_{\kappa,m,k}
d_\kappa^{-1}d_\mu S_{\kappa^\dagger,0}
 \tF_{\mu^\dagger,\nu}\left[ \begin{array}{cc}
  \kappa & \mu^\dagger \\
  \nu^\dagger & \lambda
\end{array} \right]_{ki}^{jm}  R_{mk} \]
\end{theorem}

The formula alone does not completely prove that $S(\lambda)$ is
determined by genus zero data, since it involves $S_{\kappa^\dagger,0}$.
However, we will in fact first argue that $S(\lambda)$ for all $\lambda\in
\Lambda$ is determined by genus zero data, in particular so is
$S_{\kappa^\dagger,0} = S(0)_{\kappa^\dagger,1}^{0,1}$.

This paper is organized as follows. We present the axioms of a
modular functor in section \ref{AxiomsMF}. In section \ref{BD} we
recall the notion of basic data as defined by Kevin Walker. The
curve operators are introduced in section \ref{CO}. In section
\ref{CODT} we establish that the Dehn twist in any curve can be
expressed as a linear combination of the curve operators
associated to the curve. This rests on Proposition \ref{Cinvt}, which is
proved in section \ref{sec7}. In section \ref{From1to0} it argued
that the Dehn twist in the $\beta$ on the genus one surface in
figure \ref{figs-3} is determined via genus zero data with
respect to the factorization along the curve $\alpha$. This combined with relations
in the mapping class group of a once marked genus one surface gives
 a proof of theorem \ref{Main}. A couple of formulae involving $F$ and $S$
are derived in section \ref{sec7} followed by a proof of Proposition
\ref{Cinvt}. In section \ref{sec8} we derive the formula for the
$S(\lambda)$-matrix.

The result of this paper is used in the final paper in the series
of three papers \cite{AU1}, \cite{AU2} and \cite{AU3}. In the
first two papers we construct modular functors from Conformal
Field Theory and in the third we identifying the resulting modular
functors which underlies the Reshetikhin-Turaev TQFT \cite{RT1},
\cite{RT2} and \cite{Tu} via the Skein theory realizations of
Blanchet, Habegger, Masbaum and Vogel \cite{BHMV1},
\cite{BHMV2} and \cite{Bl1}.

\section{The axioms for a modular functor}\label{AxiomsMF}

We shall in this section give the axioms for a modular functor.
These are due to G. Segal and appeared first in \cite{Se}. We
present them here in a topological form, which is due to K. Walker
\cite{Walker}. See also \cite{Grove}. We note that similar, but
different, axioms for a modular functor are given in \cite{Tu} and
in \cite{BB}. It is however not clear if these definitions of a
modular functor is equivalent to ours.

Let us start by fixing a bit of notation. By a closed surface we
mean a smooth real two dimensional manifold. For a closed oriented
surface $\Si$ of genus $g$ we have the non-degenerate skew-symmetric
intersection pairing
\[(\cdot,\cdot) : H_1(\Si,\Z) \times H_1(\Si,\Z) \ra \Z.\]
Suppose $\Si$ is connected. In this case a Lagrangian subspace
$L\subset H_1(\Si,\Z)$ is by definition a subspace, which is
maximally isotropic with respect to the intersection pairing. - A
$\Z$-basis $(\vec \alpha, \vec \beta) = (\alpha_1,\ldots,
\alpha_g,\beta_1, \ldots \beta_g)$ for $H_1(\Si,\Z)$ is called a
symplectic basis if
\[(\alpha_i,\beta_j) = \delta_{ij}, \quad (\alpha_i,\alpha_j) = (\beta_i,\beta_j) = 0,\]
for all $i,j = 1, \ldots, g$.

If $\Si$ is not connected, then $H_1(\Si,\Z) = \oplus_i
H_1(\Si_i,\Z)$, where $\Si_i$ are the connected components of $\Si$.
By definition a Lagrangian subspace is in this paper a subspace of
the form $L = \oplus_i L_i$, where $L_i\subset H_1(\Si_i,\Z)$ is
Lagrangian. Likewise a symplectic basis for $H_1(\Si,\Z)$ is a
$\Z$-basis of the form $((\vec \alpha^i, \vec \beta^i))$, where
$(\vec \alpha^i, \vec \beta^i)$ is a symplectic basis for
$H_1(\Si_i,\Z)$.

For any real vector space $V$, we define $PV = (V-\{0\})/\R{}_+.$

\begin{definition}\label{DefPointS}
A {\em pointed surface} $(\Si,P)$ is an oriented closed surface
$\Si$ with a finite set $P\subset \Si$ of points. A pointed surface
is called {\em \stable}if the Euler characteristic of each component
of the complement of the points $P$ is negative. A pointed surface
is called {\em \fm}if each component of $\Si$ contains at least one
point from $P$.
\end{definition}

\begin{definition}\label{DefMorPointS}
A {\em morphism of pointed surfaces} $f :(\Si_1,P_1) \ra
(\Si_2,P_2)$ is an isotopy class of orientation preserving
diffeomorphisms which maps $P_1$ to $P_2$. Here the isotopy is
required not to change the induced map of
 the first order Jet at $P_1$ to the first order Jet at $P_2$.
\end{definition}

\begin{definition} \label{msurface}

A {\em marked surface\/} $ {\Sib} = (\Si, P, V, L)$ is an oriented
closed smooth surface $\Si$ with a finite subset $P \subset \Si$ of
points with projective tangent vectors $V\in \sqcup_{p \in
P}PT_{p}\Si$ and a Lagrangian subspace $L \subset H_1(\Si,\Z)$.
\end{definition}

The notions of \stable and \fm marked surfaces are
defined just like for pointed surfaces.

In the case of genus zero we omit the Lagragian subspace from the
discussion, since it in this case can only be zero subspace.

\begin{definition} \label{mmorphism}

A {\em morphism\/} $\e f : {\Sib}_1 \to {\Sib}_2$ of marked surfaces
${\Sib}_i = (\Si_i,P_i,V_i,L_i)$ is an isotopy class of orientation
preserving diffeomorphisms $f : \Si_1 \to \Si_2$ that maps
$(P_1,V_1)$ to $(P_2,V_2)$ together with an integer $s$. Hence we
write $\e f = (f,s)$.
\end{definition}

\begin{remark} {\em Any marked surface has an underlying pointed
surface, but a morphism of marked surfaces does not quit induce a
morphism of pointed surfaces, since we only require that the
isotopies preserve the induced maps on the projective tangent
spaces. }\end{remark}

\begin{remark} {\em If in the notation above, we only specify $f : \Si_1 \ra \Si_2$, then
it is assumed that the integer $s=0$. }\end{remark}

Let $\sigma$ be Wall's signature cocycle for triples of Lagrangian
subspaces of $H_1(\Si,\R{})$ (See \cite{Wall}).

\begin{definition} \label{composition}
Let $\e f_1 = (f_1,s_1) : {\Sib}_1 \to {\Sib}_2$ and $\e f_2 =
(f_2,s_2) : {\Sib}_2 \to {\Sib}_3$ be morphisms of marked surfaces
${\Sib}_i = (\Si_i,P_i,V_i,L_i)$ then the {\it composition\/} of $\e
f_1$ and $\e f_2$ is $$ \e f_2 \e f_1 = (f_2 f_1, s_2 + s_1 -
\sigma((f_2f_1)_*L_1, f_{2*}L_2,L_3)). $$
\end{definition}

With the objects being marked surfaces and the morphism and their
composition being defined as in the above definition, we have
constructed
 the category of marked surfaces.

The mapping class group $\Gamma({\Sib})$ of a marked surface ${\Sib}
= (\Si,L)$ is the group of automorphisms of ${\Sib}$. One can prove
that $\Gamma({\Sib})$ is a central extension of the mapping class
group $\Gamma(\Si)$ of the surface $\Si$ defined by the 2-cocycle $c
: \Gamma({\Sib}) \to \mathbb Z$, $c(f_1,f_2) =
\sigma((f_1f_2)_*L,f_{1*}L,L)$. One can also prove that this cocycle
is equivalent to the cocycle obtained by considering two-framings on
mapping cylinders (see \cite{At1} and \cite{A}).

Notice also that for any morphism $(f,s) : \Sib_1 \to \Sib_2$, one
can factor
\begin{eqnarray*}
(f,s) &=& \left((\id,s') : \Sib_2 \to \Sib_2\right) \circ (f,s-s')\\
&=& (f,s-s') \circ \left((\id,s') : \Sib_1 \to \Sib_1\right).
\end{eqnarray*}
In particular $(\id,s) : {\Sib} \to {\Sib}$ is $(\id,1)^s$.

\begin{definition} \label{disjunion}
The operation of {\em disjoint union of marked surfaces} is $$
(\Si_1,P_1,V_1,L_1)
 \sqcup (\Si_2,P_2,V_2,L_2) = (\Si_1 \sqcup \Si_2,P_1 \sqcup P_2,V_1\sqcup V_2,L_1 \oplus
L_2). $$

Morphisms on disjoint unions are accordingly $(f_1,s_1) \sqcup
(f_2,s_2) = (f_1 \sqcup f_2,s_1 + s_2)$.
\end{definition}

We see that disjoint union is an operation on the category of marked
surfaces.

\begin{definition}\label{or}
Let ${\Sib}$ be a marked surface. We denote by $- {\Sib}$ the marked
surface obtained from ${\Sib}$ by the {\em operation of reversal of
the orientation}. For a morphism $\e f = (f,s) : {\Sib}_1 \to
{\Sib}_2$ we let the orientation reversed morphism be given by
 $- \e f = (f,-s) : -{\Sib}_1 \to -{\Sib}_2$.
\end{definition}

We also see that orientation reversal is an operation on the
category of marked surfaces. Let us now consider glueing of marked
surfaces.

Let $(\Si, \{p_-,p_+\}\sqcup P,\{v_-,v_+\}\sqcup V,L)$ be a marked
surface, where we have selected an ordered pair of marked points
with projective tangent vectors $((p_-,v_-),(p_+,v_+))$, at which we
will perform the glueing.

Let $c : P(T_{p_-}\Si) \ra P(T_{p_+}\Si)$ be an orientation
reversing projective linear isomorphism such that $c(v_-) = v_+$.
Such a $c$ is called a {\em glueing map} for $\Si$. Let
$\tilde{\Si}$ be the oriented surface with boundary obtained from
$\Si$ by blowing up $p_-$ and $p_+$, i.e.
\[\tilde{\Si} = (\Si -\{p_-,p_+\})\sqcup P(T_{p_-}\Si)\sqcup P(T_{p_+}\Si),\]
with the natural smooth structure induced from $\Si$. Let now
$\Si_c$ be the closed oriented surface obtained from $\tilde{\Si}$
by using $c$ to glue the boundary components of $\tilde{\Si}$. We
call $\Si_c$ the glueing of $\Si$ at the ordered pair
$((p_-,v_-),(p_+,v_+))$ with respect to $c$.

Let now $\Si'$ be the topological space obtained from $\Si$ by
identifying $p_-$ and $p_+$. We then have natural continuous maps $q
: \Si_c \ra \Si'$ and $n : \Si \ra \Si'$. On the first homology
group $n$ induces an injection and $q$ a surjection, so we can
define a Lagrangian subspace $L_c \subset H_1(\Si_c,\Z)$ by $L_c =
q_*^{-1}(n_*(L))$. We note that the image of $P(T_{p_-}\Si)$ (with
the orientation induced from $\tilde{\Si}$) induces naturally an
element in $H_1(\Si_c,\Z)$ and as such it is contained in $L_c$.

\begin{remark}{\em \label{remarkglue2}
If we have two glueing maps $c_i : P(T_{p_-}\Si) \ra P(T_{p_+}\Si),$
$i=1,2,$ we note that there is a diffeomorphism $f$ of $\Si$
inducing the identity on $(p_-,v_-)\sqcup(p_+,v_+)\sqcup(P,V)$ which
is isotopic to the identity among such maps, such that
$(df_{p_+})^{-1} c_2 df_{p_-} = c_1$. In particular $f$ induces a
diffeomorphism $f : \Si_{c_1} \ra \Si_{c_2}$ compatible with $f :
\Si \ra \Si$, which maps $L_{c_1}$ to $L_{c_2}$. Any two such
diffeomorphims of $\Si$ induces isotopic diffeomorphims from $\Si_1$
to $\Si_2$.}\end{remark}

\begin{definition} \label{glueing}
Let ${\Sib} = (\Si, \{p_-,p_+\}\sqcup P,\{v_-,v_+\}\sqcup V,L)$ be a
marked surface. Let $$c : P(T_{p_-}\Si) \ra P(T_{p_+}\Si)$$ be a
glueing map and $\Si_c$ the glueing of $\Si$ at the ordered pair
$((p_-,v_-),(p_+,v_+))$ with respect to $c$. Let $L_c \subset
H_1(\Si_c,\Z)$ be the Lagrangian subspace constructed above from
$L$. Then the marked surface ${\Sib}_c = (\Si_c,P,V,L_c)$ is defined
to be the {\em glueing} of ${\Sib}$ at the ordered pair
$((p_-,v_-),(p_+,v_+))$ with respect to $c$.
\end{definition}

We observe that glueing also extends to morphisms of marked surfaces
which preserves the ordered pair $((p_-,v_-),(p_+,v_+))$, by using
glueing maps which are compatible with the morphism in question.

\begin{remark}
Let $\Sib = (\Si, P, V, L)$ be marked surface. Assume that $\gamma$ is an oriented
simple closed curve on $\Si-P$, such that $[\gamma] \in L$. Assume
further we have a point $p$ on $\gamma$. We can then cut $\Si$
along $\gamma$ and obtain a surface with two boundary components,
which are naturally identified with $\gamma$.
By identifying each of the two boundary component to a point, say $\{p', p''\}$, we get a new
closed surface $\tSi$, with a set of marked points $\tP = O \cup \{p',
p''\}$ and tangent vectors $\tV = V \cup \{v',v''\}$. Here $v'$
and $v''$ are induced by $p\in \gamma$. Let $\Si'$ be obtained
from $\tSi$ by identifying $p'$ with $p''$. We have the quotient
map $q : \Si \ra \Si'$ and the identification map $n : \tSi \ra
\Si'$. We specify a Lagragian subspace in $\tL\subset H_1(\tSi,\bZ)$
by $\tL = n_*^{-1}q_*(L)$.

We say that $\tSib = (\tSi, \tP, \tV, \tL)$ is obtained from
$\Sib$ by factoring $\Sib$ along $(\gamma, p)$. The operation of
factoring is an inverse to glueing.
\end{remark}

We can now give the axioms for a 2 dimensional modular functor.

\begin{definition} \label{DefLS}

A {\em label set\/} $\L$ is a finite set furnished with an
involution $\l \mapsto  \l^\dagger$ and a trivial element $0$ such that
$ 0^\dagger = 0$.
\end{definition}

\begin{definition} \label{lmsurface}

Let $\L$ be a label set. The category of {\em $\L$-labeled marked
surfaces\/} consists of marked surfaces with an element of $\L$
assigned to each of the marked point and morphisms of labeled marked
surfaces are required to preserve the labelings. An assignment of
elements of $\L$ to the marked points of ${\Sib}$ is called a
labeling of ${\Sib}$ and we denote the labeled marked surface by
$({\Sib},\l)$, where $\l$ is the labeling.
\end{definition}

We define a labeled pointed surface similarly.

\begin{remark}{\em
The operation of disjoint union clearly extends to labeled marked
surfaces. When we extend the operation of orientation reversal to
labeled marked surfaces, we also apply the involution $ (\cdot)^\dagger$
to all the labels. }\end{remark}

\begin{definition} \label{DefMF}
A {\em modular functor\/} based on the label set $\L$ is a functor
$V$ from the category of labeled marked surfaces to the category of
finite dimensional complex vector spaces satisfying the axioms MF1
to MF5 below.
\end{definition}

\subsubsection*{MF1} {\it Disjoint union axiom\/}: The operation of disjoint
union of labeled marked surfaces is taken to the operation of tensor
product, i.e. for any pair of labeled marked surfaces there is an
isomorphism $$ V(({\Sib}_1,\l_1) \sqcup ({\Sib}_2,\l_2)) ) \cong
V({\Sib}_1,\l_1) \otimes V({\Sib}_2,\l_2). $$ The identification is
associative.

\subsubsection*{MF2} {\it Glueing axiom\/}: Let ${\Sib} $ and ${\Sib}_c$ be
marked surfaces such that ${\Sib}_c$ is obtained from ${\Sib} $ by
glueing at an ordered pair of points and projective tangent vectors
with respect to a glueing map $c$. Then there is an isomorphism $$
V({\Sib}_c,\lambda) \cong \bigoplus_{\m \in \L} V({\Sib},\m,
\m^\dagger,\l), $$ which is associative, compatible with glueing of
morphisms, disjoint unions and it is independent of the choice of
the glueing map in the obvious way (see remark \ref{remarkglue2}).

\subsubsection*{MF3} {\it Empty surface axiom\/}: Let $\emptyset$ denote
the empty labeled marked surface. Then $$ \dim V(\emptyset) = 1.
$$

\subsubsection*{MF4} {\it Once punctured sphere axiom\/}: Let $\Sib = (S^2,
\{p\},\{v\},0)$ be a marked sphere with one marked point. Then $$
\dim V(\Sib,\l) = \left\{ \begin{array}{ll} 1,\qquad &\l = 0\\
0,\qquad & \l \ne 0.\end{array}\right. $$

\subsubsection*{MF5} {\it Twice punctured sphere axiom\/}: Let $\Sib = (S^2,
\{p_1,p_2\},\{v_1,v_2\},\{0\})$ be a marked sphere with two marked
points. Then $$ \dim V(\Sib,(\l,\mu)) = \left\{ \begin{array}{ll} 1,
\qquad &\l =  \mu^\dagger\\ 0,\qquad &\l \ne  \mu^\dagger.\end{array}\right.
$$

In addition to the above axioms one may has extra properties, namely

\subsubsection*{MF-D} {\it Orientation reversal axiom\/}:
The operation of orientation reversal of labeled marked surfaces is
taken to the operation of taking the dual vector space, i.e for any
labeled marked surface $({\Sib},\l)$ there is a pairing $$ \langle
\cdot,\cdot\rangle : V({\Sib},\l) \otimes V(-{\Sib}, \l^\dagger) \ra
\C{}, $$ compatible with disjoint unions, glueings and orientation
reversals (in the sense that the induced isomorphisms $ V({\Sib},\l)
\cong V(-{\Sib},\dagger \l)^*$ and $V(-{\Sib},\l^\dagger) \cong
V({\Sib},\l)^*$ are adjoints).

\vskip.4cm

 and

\subsubsection*{MF-U} {\it Unitarity axiom\/}

Every vector space $V({\Sib},\l)$ is furnished with a hermitian
inner product
$$ ( \cdot,\cdot ) : V({\Sib},\l) \otimes \overline{V({\Sib},\l)}
\to {\mathbb C} $$ so that morphisms induces unitary transformation.
The hermitian structure must be compatible with disjoint union and
glueing. If we have the orientation reversal property, then
compatibility with the unitary structure means that we have a
commutative diagram $$\begin{CD} V({\Sib},\l) @>>\cong>
V(-{\Sib},\dagger \l)^*\\ @VV\cong V @V\cong VV\\
\overline{V({\Sib},\l)^*} @>\cong>> \overline{V(-{\Sib},\l^\dagger)},
\end{CD}$$
where the vertical identifications come from the hermitian structure
and the horizontal ones from the duality.

\section{Basic data}\label{BD}

Following Walker \cite{Walker}, we will review the notion of basic
data.

Fix throughout the rest of this paper a modular functor $Z$ with label
set $(\Lambda,0,\dagger)$.
We do not assume that this modular functor $Z$ has a duality structure as in axiom MF-D or a
unitary structure as in axiom MF-U.

We let $\Delta = (S^2;\infty; v_\infty)$, $\Xi = (S^2; 0,\infty;v_0, v_\infty)$
and recall that $\Upsilon = (S^2;0,1,\infty; v_0,v_1,v_\infty)$. For
$\lambda,\mu,\nu\in \Lambda$ we define
\[Z_0 = Z(\Delta, 0), \mbox{ } Z_{\lambda,\lambda^\dagger} = Z(\Xi,\lambda,\lambda^\dagger) \]
and recall that
\[Z_{\lambda,\mu,\nu} = Z(\Upsilon, \lambda,\mu,\nu).\]
Here $\dim Z_0 = \dim Z_{\lambda,\lambda^\dagger} = 1$ and we
define
\[N_{\lambda,\mu}^\nu = \dim Z_{\lambda,\mu,\nu^\dagger}.\]
The morphisms of the marked surfaces $\Delta$, $\Xi$ and $\Upsilon$
acts on the vector spaces $Z_0$, $Z_{\lambda,\lambda^\dagger}$ and
$Z_{\lambda,\mu,\nu}$.

Recall the self morphisms $R$ and $B$ of $\Upsilon$. They,
together with the Dehn-twist around $0$ and the morphisms
$(\Id,1)$, generate the mapping class groupoids of $\Upsilon$
with all its possible labelings. By factorization, these also
determine the action of the mapping class groupoid of both
$\Delta$ and $\Xi$ with all possible labelings. We shall make use of the
following notation:
$$
B_{12} = RBR^{-1}, \mbox{ } B_{13} = R^{-1}BR \mbox{ and } B_{23} = B$$
Further we denote by $T_1$, $T_2$ and $T_3$ the Dehn-twists around $0$,
$1$ and $\infty$ respectively.

Let $\Kappa = (S^2; -2,-1,1,2; v_{-2},v_{-1},v_1,v_2)$. Let
$\gamma_1$ be the circle of radius $3/2$ centered in $0$ and
$\gamma_2$ be the circle of radius $1$ centered in $3/2$, both
considered as simple closed curves on $\Kappa$. See figure \ref{figs-31}.

Factoring $\Kappa$ along $(\gamma_1,3/2)$ we obtain the marked surface
$\Kappa_1$ say and likewise factoring $\Kappa$ along $(\gamma_2,5/2)$ we
obtain the marked surface
$\Kappa_2$ say. There is a unique diffeomorphism from $\Upsilon \coprod \Upsilon$
to $\Kappa_i$ which takes the real axes to parts of the
real axes on each component, infinity to the respective quotients
of $\gamma_i$ and maps the first copy of $\Upsilon$ to the interior of $\gamma_i$ in $\C{}$.

We thus get two isomorphisms
\[\Phi_1(\mu,\lambda) : \bigoplus_{\nu \in \Lambda} Z(\Upsilon; \nu,  \mu, \lambda)
\otimes Z(\Upsilon; \nu^\dagger, \kappa, \xi ) \ra
Z(\Kappa, \lambda, \mu, \xi, \kappa)
\]
and
\[\Phi_2(\mu,\lambda) : \bigoplus_{\tnu \in \Lambda} Z(\Upsilon; \tnu,  \lambda, \kappa )
\otimes Z(\Upsilon;  \tnu^\dagger,  \xi,\mu) \ra
Z(\Kappa, \lambda, \mu, \xi, \kappa).
\]

We then define
\begin{equation}
F \left[ \begin{array}{cc}
  \mu & \xi \\
  \lambda & \kappa
\end{array} \right]: \bigoplus_{\nu \in \Lambda} Z_{\nu, \mu, \lambda  }
\otimes Z_{\nu^\dagger, \kappa, \xi } \ra
\bigoplus_{\tnu \in \Lambda}Z_{\tnu, \lambda, \kappa }
\otimes Z_{\tnu^\dagger, \xi, \mu }\label{F}
\end{equation}
by
\[F = \Phi_2^{-1}\circ \Phi_1.\]

Let $\Theta$ be an oriented genus one surface. Let $p$ be a point on $\Theta$
and $v_p$ be a tangent
vector at $p$. Choose two oriented simple closed
curves $(\alpha,\beta)$ on $\Theta-\{p\}$ as indicated in figure \ref{figs-3}.

Let $L_\alpha= \spane \{[\alpha]\}$ and $L_\beta = \spane
\{[\beta]\}$ be the Lagrangian subspaces generated in the first homology group.
Let $\Thb^\alpha =  (\Theta,p,v_p,L_\alpha)$ and $\Thb^\beta = (\Theta,p,v_p,L_\beta)$

Let $\Thb_\alpha$ and $\Thb_\beta$ be marked surfaces, which
results from factoring $\Thb^\alpha$ along $\alpha$, respectively
$\Thb^\beta$ along $\beta$.
By factorization we get isomorphisms
\[\Phi_\alpha : Z(\Thb^\alpha,\lambda) \ra \bigoplus_{\mu} Z(\Thb_\alpha,\lambda, \mu,\mu^\dagger)\]
and
\[\Phi_\beta : Z(\Thb^\beta,\lambda) \ra \bigoplus_{\nu} Z(\Thb_\beta,\lambda,\nu,\nu^\dagger).\]

Pick diffeomorphisms
$$ f_\alpha : \Upsilon \ra \Thb_\alpha$$
and
$$ f_\beta : \Upsilon \ra \Thb_\beta$$
which maps the real axis onto $\beta$, respectively $\alpha$.

Then we can define
\begin{equation}
S(\lambda) : \bigoplus_{\mu} Z_{\lambda,\mu,\mu^\dagger}
\ra \bigoplus_{\nu} Z_{\lambda,\nu,\nu^\dagger}\label{S}
\end{equation}
by
$$
S(\lambda) = Z(f_\beta)^{-1} \Phi_\beta Z(\Id) \Phi_\alpha^{-1}Z(f_\alpha),
$$
where $\Id : (\Thb^\alpha,\lambda) \ra (\Thb^\beta,\lambda)$.

Following Walker ((3.6) in \cite{Walker}), we define basic data as
follows.

\begin{definition}[Walker]
Basic data for the modular functor $Z$ consist of the following data:

\begin{description}
\item[A] The vector
spaces $Z_0, Z_{\lambda,\lambda^\dagger}$ and
$Z_{\lambda,\mu,\nu}$ together with the induced actions of the
groupoids of morphism of marked surface acting on them.
\item[B] The linear isomorphism
\[F\left[ \begin{array}{cc}
  \mu & \xi \\
  \lambda & \kappa
\end{array} \right] : \bigoplus_{\nu \in \Lambda} Z_{\nu, \mu, \lambda  }
\otimes Z_{\nu^\dagger, \kappa, \xi } \ra
\bigoplus_{\tnu \in \Lambda}Z_{\tnu, \lambda, \kappa }
\otimes Z_{\tnu^\dagger, \xi, \mu }.\]
\item[C] The linear isomorphism
$$S(\lambda) : \bigoplus_{\mu} Z_{\lambda,\mu,\mu^\dagger} \ra \bigoplus_{\nu} Z_{\lambda,\nu,\nu^\dagger}.$$
\end{description}

\end{definition}

\begin{lemma}\label{BddetZ}
The basic data determines the modular functor $Z$ uniquely.
\end{lemma}

This lemma is proved in section 5 of \cite{Walker}. Below we outline
the main construction behind that.

Of course {\bf A} and {\bf B} are by definition part of the
genus zero data of $Z$. Further {\bf A} and {\bf B} clearly determines the genus
zero part of $Z$. By definition $S(\lambda)$ requires genus one as
well. But the main result of this paper is that {\bf C} is in
fact determined by {\bf A} and {\bf B}. Hence

\begin{theorem}\label{12BddetZ}
The basic data under {\bf A} and {\bf B} determines
the modular functor $Z$.
\end{theorem}

This theorem follows from Theorem \ref{Main}, which is proved
below.

Let us now fix a vector $\zeta_0 \in Z_0$ and
$\zeta(\lambda) \in Z_{\lambda,\lambda^\dagger}$. By the glueing
axiom we get natural isomorphisms
\[Z_{0} \cong Z_{0,0} \otimes Z_0\]
and
\[Z_{\lambda,\lambda^\dagger} \cong Z_{\lambda,\lambda^\dagger} \otimes Z_{\lambda,\lambda^\dagger}.\]
Under these isomorphisms, we require that
\[\zeta_0 = \zeta(0) \otimes \zeta_0\]
and
\[\zeta(\lambda) = \zeta(\lambda) \otimes \zeta(\lambda).\]

This condition uniquely fixes $\zeta(\lambda)$ for all $\lambda \in
\Lambda$.

We shall also need to do a computation in a basis of
$Z_{\lambda,\mu,\nu}$, so we fix a basis $\zeta_{j}(\lambda,\mu,\nu)$,
$j = 1, \ldots, N_{\lambda,\mu}^{\nu^\dagger}$ for each of these vector spaces.
In case $N_{\lambda,\mu}^{\nu^\dagger} = 1$, we use the notation $\zeta(\lambda,\mu,\nu) =
\zeta_{1}(\lambda,\mu,\nu)$. We will require
that under the isomorphism
\[Z_{\lambda,\lambda^\dagger} \cong Z_{\lambda,\lambda^\dagger,0}\otimes Z_0\]
we have that
\begin{equation}
\zeta(\lambda) = \zeta_1(\lambda,\lambda^\dagger,0) \otimes
\zeta_0.\label{ADP}
\end{equation}

We write $S = S(0)$ and with respect to the prefereed basis $\zeta(\lambda)
\in Z_{\lambda,\lambda^\dagger}$, $\lambda\in \Lambda$, we have
the matrix presentation $S = (S_{\lambda,\mu})_{\lambda,\mu\in
\Lambda}$.

Let us now recall the reconstruction of a modular functor from its basic data. On a
marked surface $\Sib$, one considers pairs $(C,\Pi)$, where
\begin{itemize}
\item $C$ is a finite collections of disjoint
simple closed curves, each equipped with a base point, such the result of
factoring $\Sib$ along $C$ results in a disjoint union of marked
surfaces $\Sib_i$, $i\in I$ for some finite set $I$. Further we
assume that the Lagrangian subspace of $\Sib$ is generated by the curves in
$C$.
\item $\Pi$ is a disjoint union of
morphisms of marked surfaces from an ordered disjoint union of
$\Delta$'s, $\Xi$'s and $\Upsilon$'s to the $\Sib_i$, $i\in I$, covering each $\Sib_i$
exactly once.
\end{itemize}
Such a pair $(C,\Pi)$ is called an overmarking of $\Sib$ in
\cite{Walker}. If $C$ is such that $|C|$ is minimal, then $(C,\Pi)$ is
called a marking in \cite{Walker} following \cite{HT}. We shall call these pairs {\em
decompositions} in this paper so as not to confuse them with the
structure of a marked surface as introduced in the previous section.

It follows from the results of \cite{HT}, that any two pairs of
decompositions of a given marked surface are related by a finite
sequence of one of the following changes of decompositions from say
$(C',\Pi')$ to $(C'',\Pi'')$:
\begin{description}
\item[$\M$] The collection of curves are the same $C' = C''$, and there are
automorphisms of the $\Delta$'s, $\Xi$'s
and $\Upsilon$'s, which relates $\Pi'$ to $\Pi''$, namely $(\Pi'')^{-1}\Pi'$.
\item[$\A$] Insertion or removal of a component of $C'$ to obtain $C''$,
 which results in a corresponding
insertion or removal of a copy of $\Xi$.
\item[$\cD$] Insertion or removal of a component of $C'$ to obtain $C''$, which results
in the replacement of one $\Xi$ by one new $\Delta$ and one new $\Upsilon$ or the converse.
\item[$\cF$] There are $\gamma' \in C'$ and $\gamma''\in C''$ with the property that
 $C'- \{\gamma'\} = C'' - \{\gamma''\}$ and the factorization
 along $C'- \{\gamma'\}$ contains a component which via $\Pi'$ and
 $\Pi''$ is identified with $\Psi$, such that $\gamma_1$ goes to
 $\gamma'$ and $\gamma_2$ goes to $\gamma''$.
\item[$\cS$] There are $\gamma' \in C'$ and $\gamma''\in C''$ with the property that
 $C'- \{\gamma'\} = C'' - \{\gamma''\}$ and the factorization
 along $C'- \{\gamma'\}$ contains a component which via $\Pi'$ and
 $\Pi''$ is identified with $(\Theta,p,v_p)$, such that $\alpha$ goes to
 $\gamma'$ and $\beta$ goes to $\gamma''$.
\end{description}

Let $(C,\Pi)$ be a decomposition of a labeled marked surface $(\Sib,
\lambda)$. Let $\Sib_C$ be the marked surface, one obtains from
factoring $\Sib$ along $C$. Let $\Lambda_c = \Lambda^{\times c}$, where
$c = |C|$. The factorization axiom gives us
an isomorphism
\[Z(\Sigma,\lambda) \cong \bigoplus_{\mu \in \Lambda_{C}} Z(\Sigma_C,\lambda,\mu)\]
For each $\mu \in \Lambda_c$, we let $Z(\lambda,\mu)$ be the
corresponding tensor product of the vector spaces $Z_{\lambda'}$'s,
$Z_{\lambda',(\lambda')^\dagger}$'s and
$Z_{\lambda',\mu',\nu'}$'s.
Using $\Pi$ we then get an isomorphism
$$Z(\Pi) : Z(\lambda,\mu) \ra Z(\Sigma_C,\lambda,\mu)
$$
which induces an isomorphism
$$
Z(C,\Pi) : \bigoplus_{\mu \in \Lambda_{c}} Z(\lambda,\mu) \ra Z(\Sigma,\lambda).$$

Suppose now $(C',\Pi')$ and $(C'',\Pi'')$ are two decompositions of
the same pointed surface $(\Si,P,V)$, which are
related by one of the changes $\M$ to $\cS$. Let
$L'$ and $L''$ be the Lagrangian subspaces generated by respectively the
$C'$ and $C''$. In the cases $\M$ to $\cF$, we have that $L' = L''$, which is
not the case for $\cS$. Let $\Sib' = (\Si,P,V,L')$ and $\Sib'' =
(\Si,P,V,L'')$. We then get an induced isomorphism
$$Z((C',\Pi'), (C'',\Pi'')) : \bigoplus_{\mu \in \Lambda_{c'}}
Z(\lambda,\mu)\ra \bigoplus_{\mu \in \Lambda_{c''}}
Z(\lambda,\mu)$$
given by
$$
Z((C',\Pi'), (C'',\Pi''))  =  Z(C'',\Pi'')^{-1}Z(\Id: (\Sib',\lambda) \ra
(\Sib'',\lambda)) Z(C',\Pi')$$
This isomorphism is determined by the basic data as
follows:

For change of type:

\begin{description}
\item[$\M$] The linear map is  given by a direct sum of tensor products
of the linear maps induced by the morphism $(\Pi'')^{-1} \circ \Pi'$ between
the appropriate $\Xi$'s, $\Delta$'s and $\Upsilon$'s.
\item[$\A$] The linear map is induced by insertion or "removal" of the
vector $\zeta(\lambda) \in Z_{\lambda,\lambda^\dagger}$ for the
appropriate $\Xi$.
\item[$\cD$] The linear map is induced by the identity tensor the isomorphism
(\ref{ADP}) (or its inverse) inserted at the appropriate place.
\item[$\cF$] The linear map is induced by the identity tensor the isomorphism
(\ref{F}) inserted at the appropriate place.
\item[$\cS$] The linear map is induced by the identity tensor the isomorphism
(\ref{S}) inserted at the appropriate place.
\end{description}

Let $(C_1,\Pi_1)$ and $(C_2,\Pi_2)$ be any two decompositions of the same labeled pointed
surface $(\Si, P, V,\lambda)$. Let $L_i$ be the induced Lagrangian subspaces on $\Si$
and $\Sib_i = (\Si, P, V, L_i)$ the corresponding marked surfaces. Since
$(C_1,\Pi_1)$ and $(C_2,\Pi_2)$ are
related by a sequence of changes $\M$ to $\cS$, we get that the
isomorphism
$$Z((C_1,\Pi_1), (C_2,\Pi_2)) : \bigoplus_{\mu_1 \in \Lambda_{c_1}}
Z(\lambda,\mu_1)\ra \bigoplus_{\mu_2 \in \Lambda_{c_2}}
Z(\lambda,\mu_2)$$
given by
$$
Z((C_1,\Pi_1), (C_2,\Pi_2))  =  Z(C_2,\Pi_2)^{-1}Z(\Id: (\Sib_1,\lambda) \ra
(\Sib_2,\lambda)) Z(C_1,\Pi_1)$$
is also determined by the basic data.

Suppose now that $f : (\Sib_1,\lambda_1) \ra (\Sib_2,\lambda_2)$ is a morphism of marked
surfaces and $(C,\Pi)$ is a decomposition of $\Sib_1$. Then $(f(C), f\circ
\Pi)$ is a decomposition of $\Sib_2$ and we get a commutative
diagram:
\[
\xymatrix{
  \bigoplus_{\mu \in \Lambda_{c}} Z(\lambda_1,\mu) \ar[d]_{\Id} \ar[r]^{Z(C,\Pi)} & Z(\Sib_1,\lambda_1)
   \ar[d]^{Z(f)} \\
  \bigoplus_{\mu \in \Lambda_{c}} Z(\lambda_2,\mu) \ar[r]^{Z(f(C), f\circ \Pi)} & Z(\Sib_2,\lambda_2)
  }.
\]

Suppose now $(\gamma, p)$ is an oriented simple closed curve with
a preferred point $p$ on a marked surface $\Sib$. Let $\tSib$ be
obtained from $\Sib$ by factoring along $(\gamma, p)$. Suppose
that $(C,\Pi)$ is a decomposition of $\Sib$, such that $\gamma \in
C$. Then $(C,\Pi)$ also induces a decomposition of $\tSib$, say $(\tC, \tPi)$
and we get the following commutative diagram
\[
\xymatrix{
  \bigoplus_{\mu \in \Lambda_{c}} Z(\lambda,\mu) \ar[d]_{\Id} \ar[r]^{Z(C,\Pi)} & Z(\Sib,\lambda)
   \ar[d]^{\cong} \\
  \bigoplus_{\mu' \in \Lambda} \bigoplus_{\mu'' \in \Lambda_{c-1}}
  Z(\lambda,\mu',\mu'') \ar[r]^{Z(\tC, \tPi)} &
  \bigoplus_{\mu' \in \Lambda}Z(\tSib,\lambda,\mu',(\mu')^\dagger)
  }.
\]
Likewise we trivially get a similar diagram for the case of the
disjoint union isomorphism.

From this it follows that the basic data determines the modular
functor because the action of any morphism, the factorization
and the disjoint union isomorphisms are all determined by the
basic data.

By considering various finite sequences of decompositions of certain
marked surfaces, we generate nontrivial relations on the basic
data. Ten such universal relations are given on page 55 in \cite{Walker}.
We will use some of them in section \ref{sec7}.

\section{Curve operators}\label{CO}

Let $\Sib = (\Si, P, V, L)$ be a general marked surface. Let
$\lambda$ be a labeling of $\Sib$.

Let $\gamma$ be an oriented simple closed curve on $\Sigma - P$ and $\lg
\in \Lambda$ a fixed label. We now define an operator
\[Z(\ga, \lg) : Z(\Sib,\lambda) \ra Z(\Sib,\lambda)\]
canonically associated to the pair $(\ga,\lg)$.

\begin{center}
\begin{figure}[H]
\includegraphics[scale=.5]{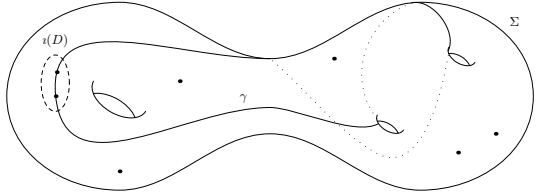}
\caption{The curve $\gamma$ and the disk $\imath(D)$ on $\Si$.}
\label{figs-4}
\end{figure}
\end{center}

Choose an embedding $\imath : D \ra \Si-P$ of the unit disc $D$ into
$\Si-P$, such that $\imath([-1,1]) = \gamma\cap \imath(D)$ as indicated
in figure \ref{figs-4}.
Let $(\ga_\imath,p_\imath)
= (\imath(\partial D), \imath(1))$. Let $P' = \{\imath(-\frac12), \imath(\frac12)\}$ and
$V'$ be directions along $\gamma$ in the positive direction at
$P'$. Let $\tP = P \cup P'$, $\tV = V \cup V'$ and $\tSib = (\Si,
\tP, \tV, L)$.

The factorization of  $\tSib$ along $\gamma_\imath$ has two connected
components  which we denote $\tSi'$ and $\tSi''$. Here $\tSi'$ is obtained from
$\imath(D)$ by identifying $\gamma_\imath$ to a point $p'$ with $p_\imath\in
\gamma_\imath$ inducing a tangent direction $v'$ at $p'$. Likewise
$\tSi''$ is the quotient of $\Sigma - \imath(D-\partial D)$, where we
identify $\gamma_\imath$ to a point $p''$ again with $p_\imath\in \gamma_\imath$
inducing a tangent direction $v''$ at $p''$.

Let $\tP' = P' \cup \{p'\}$ and $\tV' = V' \cup \{v'\}$ set
$\tSib' = (\tSi', \tP', \tV')$. Let $\tP'' = P \cup \{p''\}$ and
$\tV' = V \cup \{v''\}$ set $\tSib'' = (\tSi'', \tP'', \tV'')$.

The glueing and disjoint union axiom gives isomorphisms
\[Z(\tSib,\lambda, \lg^\dagger,\lg) \cong
\bigoplus_{\mu \in \Lambda} Z(\tSib', \lg^\dagger,\lg,\mu) \otimes
Z(\tSib'',\mu^\dagger, \lambda)\]
and
\[Z(\tSib,\lambda, \lg,\lg^\dagger) \cong
\bigoplus_{\mu \in \Lambda} Z(\tSib', \lg,\lg^\dagger,\mu) \otimes
Z(\tSib'',\mu^\dagger, \lambda).\]
The embedding $\imath$ induces an isomorphism of marked curves $\imath : \Upsilon \ra
\tSib'$, which therefore gives isomorphisms
\[Z(\imath) : Z(\Upsilon,\lg^\dagger,\lg, 0) \ra Z(\tSib', \lg^\dagger,\lg,0)\]
and
\[Z(\imath)^\dagger : Z(\tSib', \lg,\lg^\dagger,0)^* \ra Z(\Upsilon,\lg,\lg^\dagger, 0)^*.\]

Also by the glueing and disjoint union axiom combined with axiom
MF4, we get an isomorphism
\begin{equation}
Z(\tSib'', 0, \lambda) \cong Z(\Sib, \lambda) \label{POV}
\end{equation}
which is unique up to scale.

The vector $\zeta_1(\lg^\dagger,\lg,0)\in Z(\Upsilon, \lg^\dagger,\lg,0)$
 together
with
 the isomorphisms $Z(\imath)$ gives us
 an inclusion
\[I_\imath(\lg) :
Z(\Sib, \lambda) \ra
Z(\tSib,\lambda, \lg^\dagger,\lg) .\]
A vector vector $\alpha_\gamma \in Z(\Upsilon,
\lg,\lg^\dagger,0)^*$ together with the isomorphism
 $Z(\imath)^\dagger$ gives a projection
\[P_\imath(\lg) : Z(\tSib,\lambda, \lg,\lg^\dagger) \ra
Z(\Sib, \lambda). \]
We shall normalize the forms $\alpha_\mu\in Z(\Upsilon,
\mu,\mu^\dagger,0)^*$, $\mu\in\Lambda$ as follows:

We require that
\[\alpha_\mu(Z(B)(\zeta_1(\mu\dagger,\mu,0)))
= \frac{S_{0,\mu}}{S_{0,0}}.\]

We shall now consider the following diffeomorphism $\Phi$ of
$\tSib$. Fix a tubular neighbourhood of $\gamma$ inside $\Si - P$. The diffeomorphism
$\Phi$ will be the identity outside this tubular neighbourhood.
Inside the tubular neighbourhood it rotates and stretches
$\imath([-\frac12,\frac12])$ onto $\gamma - \imath((-\frac12,\frac12))$ and
$\gamma - \imath((-\frac12,\frac12))$ onto
$\imath([-\frac12,\frac12])$. We see that
\[Z(\Phi) : Z(\tSib,\lambda, \lg^\dagger,\lg) \ra Z(\tSib,\lambda, \lg,\lg^\dagger). \]

\begin{definition}\label{CurveOp}
The {\em curve operator} associate to $(\gamma, \lg)$ is
by definition
\[Z(\gamma, \lg) = P_\imath(\lg) \circ Z(\Phi) \circ I_\imath(\lg).\]
\end{definition}

We observe that $Z(\gamma,\lg)$ does not depend on the choice of
$\imath$. In fact
$Z(\gamma,\lg)$ only depends on the free homotopy class of
$\gamma$.

We clearly have the following lemma.

\begin{lemma}
Suppose $\e f : {\Sib}_1 \to {\Sib}_2$ is a morphism of marked
surface and $\gamma_i$ are closed oriented curves on $\Si_i -
P_i$, $i=1,2$, such that $f(\gamma_1) = \gamma_2$. Then
\[Z(\e f)^{-1} Z(\gamma_2, \lambda) Z(\e f) = Z(\gamma_1, \lambda)\]
for all $\lambda \in \Lambda$.

\end{lemma}

If $\gamma$ is contractible on the marked surface $\Si$, then
\[Z(\gamma, \lg) =  \frac{S_{0,\lg}}{S_{0,0}} \Id_{Z(\Si, \lambda)}. \]

\section{The relation between curve operators and Dehn Twists}
\label{CODT}

In this section we give a formula for the Dehn-Twist operator in
terms of the curve operators associated to any oriented simple
closed curve.

Let $(\Sib,\lambda)$ be a labeled marked surface. Let $\gamma$ be
an oriented simple closed curve on $\Si - P$. Let $\varphi_\gamma :
(\Sib,\lambda) \ra (\Sib,\lambda)$ be the Dehn twist
in the curve $\gamma$. By construction $\varphi_\gamma$ is the
identity outside some tubular neighbourhood of $\gamma$. Similarly, $Z(\gamma,
\lg)$ is also a local construction within a tubular neighbourhood
of $\gamma$. Pick a point $p$ on $\gamma$ and let $\tSib$ be obtained
from $\Sib$, by factoring $\Sib$ along $(\gamma,p)$.
Then we get an isomorphism
\[Z(\Sib,\lambda) \cong \bigoplus_{\mu}Z(\tSib,\lambda,\mu,\mu^\dagger).\]
Both of the operators $Z(\varphi_\gamma)$ and $Z(\gamma,
\lg)$ are diagonal with respect to this direct sum decomposition
and acts by multiples of the identity on each of the summands.
This follows immediately from factoring along the boundary of a
tubular neighborhood of $\gamma$. One also sees this way that these multiples by which
these operators acts by are independent of both $\gamma$ and
$(\Sib,\lambda)$.

\begin{proposition} \label{CODTp}
There exist uniquely determined constants $c_{\lg}\in \bC$ such that
for any simple closed oriented curve $\gamma$ on any labeled marked
surface $(\Sib,\lambda)$ we have that
\[Z(\varphi_\gamma) = \sum_{\lg \in \Lambda} c_{\lg} Z(\gamma, \lg).\]
\end{proposition}

Recall that $\Xi = (S^2; 0,\infty;v_0, v_\infty)$ and let $\gamma$ be
the unit circle oriented in the positive direction.

We define $C_{\lambda,\mu}\in {\mathbb C}$ to be the scalar by which
$Z(\gamma,\lambda)$ acts on $Z(\Xi,\mu,\mu^\dagger) = Z_{\mu,\mu^\dagger}$,
for $\lambda,\mu \in
\Lambda$. Further we recall that $d_\mu \in {\mathbb C}$ is the scalar by
which $Z(\varphi_\gamma)$ acts on $Z(\Xi,\mu,\mu^\dagger)= Z_{\mu,\mu^\dagger}$.

\begin{proposition}\label{Cinvt}
The matrix $C_{\lambda,\mu}$ is invertible.
\end{proposition}

We will prove this Proposition in section \ref{sec7}.

\proof[Proof of Proposition \ref{CODTp}] The constants $c_\lambda$,
which we seek has to satisfy
\[d_\mu = \sum_\lambda c_\lambda C_{\lambda,\mu}.\]
By Proposition \ref{Cinvt} the matrix $C$ is invertible, so there is
a unique solution to this set of equations. \eproof

We remark that since $C$ is invertible, the $c_\lambda$ are determined
by the genus zero data. Formula (\ref{cerne}) below gives an explicit
formula for the $c_\lambda$'s.

\section{The reduction from a once punctured genus one surfaces to genus zero}

\label{From1to0}

Recall the simple
closed curves $\alpha$ and $\beta$ on $\Theta - \{p\}$ from figure
\ref{figs-3}. Let $S :  \Thb^\alpha \ra \Thb^\alpha$ be the
morphism of marked surfaces, which satisfies that
\[S(\alpha) = \beta, \mbox{ } S(\beta) = \alpha^{-1},\]
where we here interpret $\alpha$ and $\beta$ as generators of the
fundamental group of $\Theta - \{p\}$ base at their intersection
point.

\begin{theorem}\label{Sdetbyg0}
The morphism
\[Z(S) = Z(\Thb^\alpha, \lambda) \ra Z(\Thb^\alpha, \lambda)\]
is determined by the genus zero part of $Z$ for all $\lambda\in
\Lambda$.
\end{theorem}

\proof

We recall that the mapping class group $\Gamma$ of $\Theta-\{p\}$
is
\[\Gamma \cong \{S,T | (ST)^3 = S^2\},\]
where $S$ is as specified above and $T$ is the Dehn-twist in
$\alpha$.

We see that
\[S = T^{-1}ST^{-1}S^{-1} T^{-1}.\]
Hence if we let $T'$ be the Dehn-twist in $\beta$, then
\[S = T^{-1}(T')^{-1}T^{-1}.\]
But by Proposition \ref{CODT} we have that there are constants
$\tc_\lambda$ such that
\[Z(T')^{-1} = \sum_{\lambda_\beta \in \Lambda} \tc_{\lambda_\beta} Z(\beta,\lambda_\beta). \]
So
\[Z(S) = \sum_{\lambda_\beta \in \Lambda} \tc_{\lambda_\beta}
Z(T^{-1}) Z(\beta,\lambda_\beta)Z(T^{-1}).\]

\begin{center}
\begin{figure}[H]
\includegraphics[scale=.5]{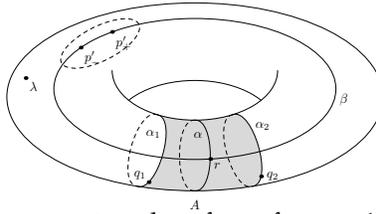}
\caption{A once punctured surface of genus 1 with an annulus A around the $\alpha$ curve.}
\label{annulus-1}
\end{figure}
\end{center}

Let $\jmath : A \ra \Theta -\{p\}$ be an embedding of an annulus
$A$ into
$\Theta-\{p\}$ as shown in figure \ref{annulus-1}.

Let
$\alpha_1 $ and $\alpha_2$ be the two boundary curves of $\jmath(A)$, with base points
say $q_1\in \alpha_1$ and $q_2\in \alpha_2$.
Let $\Theta_1$ respectively $\Theta_2$
be obtained from $\Theta$ by factoring along $\alpha_1$
respectively along $\alpha_2$. We denote the resulting base point and tangent direction by
$(P_1,V_1)=(\{q'_1,q''_1\}, \{v_{q'_1},v_{q''_1}\})$ and
$(P_2,V_2)=(\{q'_2,q''_2\}, \{v_{q'_2},v_{q''_2}\})$
on respectively $\Theta_1$ and $\Theta_2$. Let $\tSi$ be
obtained from $\Theta_1$, by factoring
along $\alpha_2$ or equivalently from $\Theta_2$ by factoring along $\alpha_1$.
Pick a point $r\in \jmath (A - \partial A) \cap
\gamma$ and the tangent direction $v_r$ at $r$ along $\beta$ in the positive direction.

Now choose an embedding $\imath : D \ra \Theta - (\jmath(A)\cup \{p\})$ and
a diffeomorphism $\Phi$ with the properties required
in the construction of $Z(\beta,\lambda_\beta)$. In fact we will choose a
$\Phi$ which is a composite of two diffeomorphisms $\Phi_1$ and
$\Phi_2$ as follows. Let $p'_- = \imath(-\frac12)$ and $p'_+ =
\imath(\frac12)$ and  $v_{p'_-}$ respectively $v_{p'_+}$ be the induced
tangent directions along $\beta$.

Let $P = \{p, p'_-, p'_+\}$ and $V =
\{v_p, v_{p'_-}, v_{p'_+}\}$. Let $P' =  \{p, p'_+, r\}$ and $V' =
\{v_p, v_{p'_+},  v_{r}\}$. Further set
\[\tP_i = P \cup P_i, \mbox{ } \tV_i = V \cup V_i,\]
\[\tP'_i = P' \cup P_i, \mbox{ } \tV'_i = V' \cup V_i,\]
for $i=1,2$.

We pick
\[\Phi_1 : (\Theta_2, \tP_2, \tV_2) \ra (\Theta_2, \tP'_2, \tV'_2)\]
such that $\Phi_1$ is the identity outside a tubular neighbourhood of the
piece of $\beta$ from $p'_-$ to $r$ and it maps $\Phi_1(p'_-) =
p'_+$ and $\Phi_1(p'_+) =
r$. The map
\[\Phi_2 : (\Theta_1, \tP'_1, \tV'_1) \ra (\Theta_1, \tP_1, \tV_1)\]
is chosen such that it is the identity outside a neighbourhood of
the piece of $\beta$ from $r$ to $p'_-$ and it maps $\Phi_2(r) =
p'_-$. By re-glueing $\Theta_1$ and $\Theta_2$ to obtain $\Theta$, we see
that $\Phi_1$ and $\Phi_2$ induces diffeomorphisms  $\Phi'_1$ and $\Phi'_2$
of $\Theta$,
which
preserves $\alpha_2$ respectively $\alpha_1$. We let
\[\Phi = \Phi'_2 \circ \Phi'_1 : (\Theta, P, V) \ra (\Theta, P, V).\]

Thus we have the linear maps for all choices of $\mu, \mu' \in \Lambda$
\[Z(\Phi_1)_\mu : Z(\Theta_2, \tP_2, \tV_2,\lambda, \lg^\dagger,\lg,\mu,\mu^\dagger)
\ra Z(\Theta_2, \tP'_2, \tV'_2,\lambda, \lg^\dagger,\lg,\mu,\mu^\dagger) \]
and
\[Z(\Phi_2)_{\mu'} : Z(\Theta_1, \tP'_1, \tV'_1,\lambda, \lg^\dagger,\lg,
\mu',(\mu')^\dagger)
\ra Z(\Theta_1, \tP_1, \tV_1,\lambda, \lg, \lg^\dagger,\mu',(\mu')^\dagger). \]

We further have the following two commutative diagrams
\[\xymatrix{
  Z(\Theta, P, V, \lambda, \lg^\dagger,\lg) \ar[d]_{} \ar[r]^{Z(\Phi'_1)} &
  Z(\Theta, P', V', \lambda, \lg^\dagger,\lg)\ar[d]^{} \\
  \bigoplus_{\mu\in \Lambda}
Z(\Theta_2,\tP_2, \tV_2,\lambda, \lg^\dagger,\lg,\mu,\mu^\dagger)
\ar[r]^{\bigoplus_{\mu\in \Lambda} Z(\Phi_1)_\mu} & \bigoplus_{\mu\in \Lambda} Z(\Theta_2,\tP'_2, \tV'_2,\lambda,
\lg^\dagger,\lg,\mu,\mu^\dagger)   }
\]

and
\[\xymatrix{
  Z(\Theta, P', V', \lambda, \lg,\lg^\dagger) \ar[d]_{} \ar[r]^{Z(\Phi'_2)} &
  Z(\Theta, P, V, \lambda, \lg,\lg^\dagger)\ar[d]^{} \\
  \bigoplus_{\mu'\in \Lambda}
Z(\Theta_1,\tP'_1, \tV'_1,\lambda, \lg,\lg^\dagger,\mu',(\mu')^\dagger)
\ar[r]^{\bigoplus_{\mu'\in \Lambda} Z(\Phi_2)_{\mu'}} &
\bigoplus_{\mu'\in \Lambda} Z(\Theta_1,\tP_1, \tV_1,\lambda, \lg,\lg^\dagger,\mu',(\mu')^\dagger)  }
\]
where the vertical arrows are the factorization isomorphisms.

By the commutativity of factorization we also have the following
commutative diagram of isomorphisms:

\[\xymatrix{
  Z(\Theta, P', V', \lambda, \lg,\lg^\dagger) \ar[d]_{} \ar[r]^{} &
 \bigoplus_{\mu\in \Lambda}
 Z(\Theta_1, \tP'_1, \tV'_1, \lambda, \lg,\lg^\dagger,\mu,\mu^\dagger)\ar[d]^{} \\
  \bigoplus_{\mu'\in \Lambda}
Z(\Theta_2,\tP'_2, \tV'_2,\lambda, \lg,\lg^\dagger,\mu',(\mu')^\dagger)
\ar@{.>}[ur]|-{} \ar[r]^{} &
 \bigoplus_{\mu, \mu'\in \Lambda}Z(\tSi,\tP, \tV,\lambda, \lg,\lg^\dagger,\mu,\mu^\dagger,\mu',(\mu')^\dagger)   }
\]

Now consider the curve operator $Z(\beta,\lambda_\beta)$. Since
factorization along non-intersecting curves commute we see that
both  $P_\imath(\lambda_\beta)$ and $I_\imath(\lambda_\beta)$ is
determined by genus zero morphism by factoring along say
$\alpha_1$, which commutes with the factorization along $\gamma_\imath$.

By representing $T$ as the Dehn Twist in $\alpha_1$ respectively
$\alpha_2$, we see get diffeomorphisms
\[T_1 : (\Theta_2, \tP_2',\tV_2') \ra (\Theta_2, \tP_2',\tV_2')\]
and
\[T_2 : (\Theta_1, \tP_1',\tV_1') \ra (\Theta_1, \tP_1',\tV_1').\]
Since both of these diffeomorphisms are the identity in a
neighbourhood of $\imath(D)$, we get the formula
\[Z(S) = \sum_{\lambda_\beta \in \Lambda} \tc_{\lambda_\beta}
 P_\imath(\lambda_\beta)
 Z(T_2^{-1}\Phi'_2 ) Z(\Phi'_1 T_1^{-1}) I_\imath(\lambda_\beta).\]
Tracing through the previous three commutative
diagrams, we see that this is also the case for $Z(\Phi_1'T_1^{-1})$ and
$Z(T_2^{-1}\Phi_2')$ are determined by genus zero data.
\eproof

We observe the same argument can be used to show the following:
The action of a Dehn twist along any curve on a marked surface equipped
with a fixed decomposition is determined by the isomorphisms $\M$ to $\cF$.

\proof[Proof of Theorem \ref{Main}]
We consider the morphism
\[S : (\Theta^\alpha, \lambda) \ra (\Theta^\beta, \lambda)\]
and observe that it is compatible with the factorizations in
$\alpha$ and $\beta$ and via $f_\alpha$ and $f_\beta$ is
compatible with the identity morphism $\Id : \Upsilon \ra
\Upsilon$. We have the composition identity
\begin{eqnarray*}
\lefteqn{S : (\Theta^\alpha, \lambda) \ra (\Theta^\beta,
\lambda)}\\
&  = &
 (S : (\Theta^\alpha, \lambda) \ra (\Theta^\beta, \lambda)) \circ
 (\Id : (\Theta^\beta, \lambda) \ra (\Theta^\alpha, \lambda)).
\end{eqnarray*}

From this we conclude that
\[Z(f_\alpha^{-1}) Z(S: (\Theta^\alpha, \lambda) \ra (\Theta^\alpha, \lambda)) Z(f_\alpha)
= S(\lambda)^{-1}.\]
Since $Z(S) :  Z(\Theta^\alpha, \lambda) \ra Z(\Theta^\alpha,
\lambda)$ is determined by genus zero data, we see that so is
$S(\lambda)$ for all $\lambda \in \Lambda$.
\eproof

\section{Consequences of the universal relations on the basic data.}
\label{sec7}

In this section we will prove the formula in Proposition
\ref{Cinvt}. It follows immediately from Proposition
\ref{CurveOpformula} below. However, first we need to derive a
couple of
consequences from the pentagon relation for $F$ and
relations between $F$ and $S$.

Consider the isomorphism (\ref{F}). A special case is if one of the labels $\lambda,\mu, \xi,
\kappa$ equals $0\in \Lambda$. In this case all of the terms in the direct
 sum of the domain and codomain of (\ref{F}) are zero, except for one of them.

\begin{lemma}
We have the formulae
\begin{description}
\item[$\lambda = 0$]  $F(\zeta(\mu^\dagger,0,\mu) \otimes v) =
R^2(v) \otimes \zeta(\xi^\dagger,\xi,0)$
\item[$\mu = 0$]  $F(\zeta(\lambda^\dagger,\lambda,0) \otimes v) =
\zeta(\kappa^\dagger,0,\kappa) \otimes R(v)$
\item[$\xi = 0$]  $F(v \otimes \zeta(\kappa^\dagger,\kappa,0)) =
R(v) \otimes \zeta(\lambda^\dagger,0,\lambda)$
\item[$\kappa = 0$]  $F(v \otimes \zeta(\xi^\dagger,0,\xi)) =
\zeta(\mu^\dagger,\mu,0) \otimes R^2(v)$
\end{description}
\end{lemma}

\proof
The case $\xi=0$ is precisely relation 5. on page 55 of \cite{Walker}. They are all
obtained by considering the two ways of decomposing $\Upsilon$
along a system of two curves, related by an $\cF$ change and such that the factorization
of both gives two copies of $\Upsilon$ and one copy of $\Delta$, as illustrated in
figure 27 on page 53 of \cite{Walker}.
\eproof

For each $\lambda \in \Lambda$ we define
\[E_\lambda = F_{0,0}\left[ \begin{array}{cc}
  \lambda & \lambda^\dagger \\
  \lambda^\dagger & \lambda
\end{array} \right].\]

\begin{lemma}\label{ESS}
We have the following formula
\[S_{0,0} E_\lambda = S_{0,\lambda^\dagger}.\]
for all $\lambda \in \Lambda$.
\end{lemma}

\begin{center}
\begin{figure}[H]
\includegraphics[scale=.5]{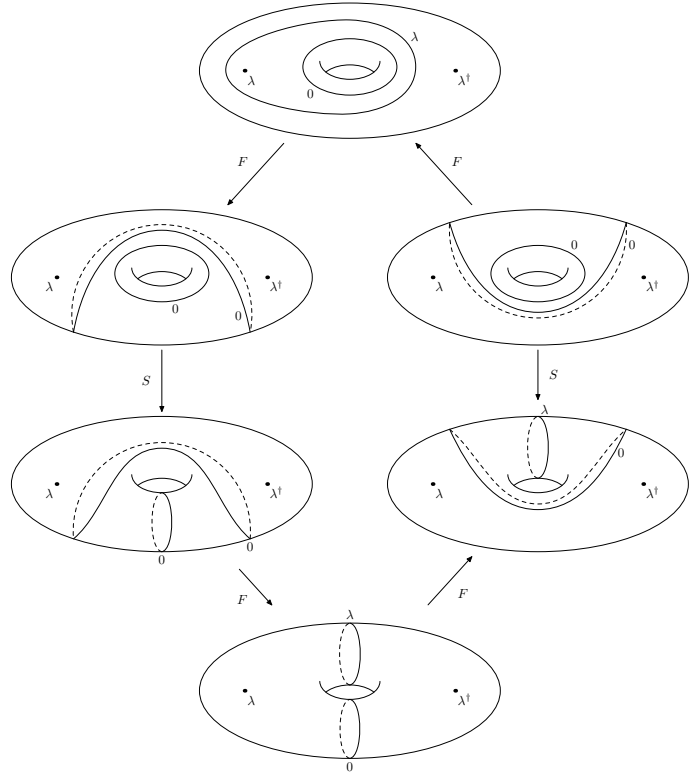}
\caption{A relation between $S$ and $F$.}
\label{SF}
\end{figure}
\end{center}

\proof We consider relation 3. on page 55 of \cite{Walker},
which is obtained from considering six decompositions of a
genus one surface with two marked points as depicted in figure
\ref{SF}. See also figure 18 on page 43 and figure 26 on page 52 of
\cite{Walker}.

Starting at the domain of the morphism $S$ on the right of figure \ref{SF} and
going counter clockwise around to the codomain of the same
morphism, we get that
\begin{eqnarray*}
\lefteqn{\zeta(0,\lambda^\dagger,\lambda) \otimes \zeta(0,0,0)
\stackrel{F}{\longmapsto}}\\
& & \zeta(\lambda^\dagger,\lambda,0) \otimes
\zeta(\lambda,0,\lambda^\dagger)
\stackrel{B_{23}^{-1}\otimes B_{23}^{-1}}{\longmapsto}\\
& & \zeta(\lambda^\dagger,0,\lambda) \otimes
\zeta(\lambda,\lambda^\dagger,0)
\stackrel{F}{\longmapsto}\\
& & \zeta(0,\lambda,\lambda^\dagger) \otimes \zeta(0,0,0)
\stackrel{\Id\otimes T_{3}^{-1}T_1 B_{23}}{\longmapsto}\\
& & \zeta(0,\lambda,\lambda^\dagger) \otimes \zeta(0,0,0)
\stackrel{\Id\otimes S}{\longmapsto}\\
& & S_{0,0}\zeta(0,\lambda,\lambda^\dagger) \otimes \zeta(0,0,0)
\stackrel{F}{\longmapsto}\\
& & S_{0,0}\zeta(\lambda,\lambda^\dagger,0) \otimes \zeta(\lambda^\dagger,0,\lambda)
\stackrel{R^{-1}\otimes R}{\longmapsto}\\
& & S_{0,0}\zeta(0, \lambda,\lambda^\dagger) \otimes \zeta(0,\lambda,\lambda^\dagger)
\stackrel{F}{\longmapsto}\\
& & S_{0,0}E_\lambda\zeta(0, \lambda^\dagger,\lambda) \otimes
\zeta(0,\lambda^\dagger,\lambda).
\end{eqnarray*}
Where as the morphism $S$ of course gives
\begin{eqnarray*}
\lefteqn{\zeta(0,\lambda^\dagger,\lambda) \otimes \zeta(0,0,0)
\stackrel{\Id\otimes S}{\longmapsto}}\\
& & S_{0,\lambda^\dagger} \zeta(0,\lambda^\dagger,\lambda) \otimes
\zeta(0,\lambda^\dagger,\lambda)
\end{eqnarray*}
From which we conclude the formula
\[S_{0,0} E_\lambda = S_{0,\lambda^\dagger}.\]

\eproof

\begin{corollary}
We have that
$S_{0,0} \neq 0$
for any modular functor.
\end{corollary}

\begin{lemma}\label{Pentagonimp}
For all $\lambda,\mu,\nu \in \Lambda$, we have the relation
\begin{equation}
E_{\lambda^\dagger} \sum_{i,j,r} F_{0,\mu}\left[ \begin{array}{cc}
  \lambda & \nu \\
  \lambda^\dagger & \nu^\dagger
\end{array} \right]_{11}^{ij} R_{jr}F_{\nu,0}\left[ \begin{array}{cc}
  \lambda & \lambda^\dagger \\
  \mu^\dagger & \mu
\end{array} \right]^{11}_{rl}R^2_{is} = \delta_{sl}.\label{RelPent}
\end{equation}
\end{lemma}

\begin{center}
\begin{figure}
\includegraphics[scale=.5]{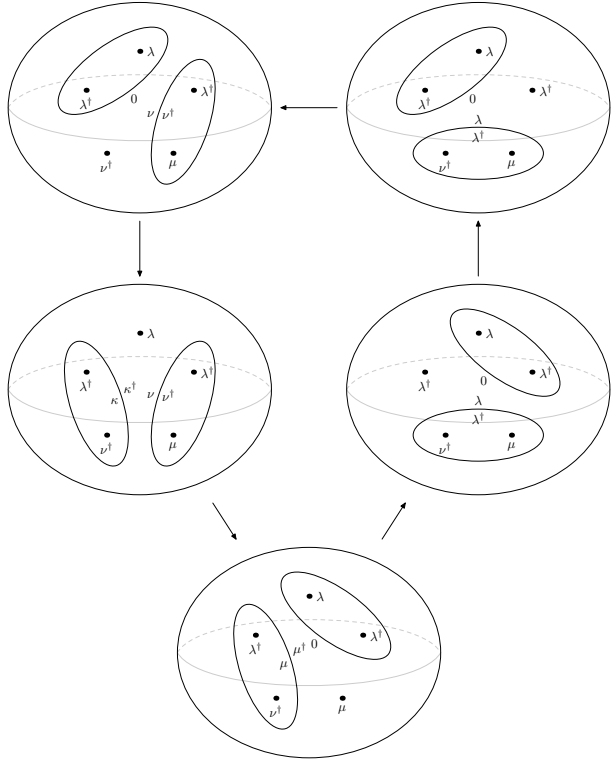}
\caption{The pentagon relation for $F$.}
\label{figs-25}
\end{figure}
\end{center}

\proof The pentagon relation as depicted in figure \ref{figs-25},
gives a relation between five applications of the $F$-isomorphism
between five decompositions of a genus zero surface
with five marked points. This is relation 1. on page 55 of \cite{Walker}
which is also illustrated in figure 24 on page 50 of \cite{Walker}.

Starting in the upper left hand corner and going counter clockwise
around we get that
\begin{eqnarray*}
\lefteqn{\zeta(0,\lambda,\lambda^\dagger)\otimes
\zeta(0,\nu,\nu^\dagger)\otimes
\zeta_l(\nu^\dagger,\mu,\lambda^\dagger) \stackrel{F\otimes \Id}{\longmapsto}} \\
& & \sum_{\kappa, i,j} F_{0,\kappa}\left[ \begin{array}{cc}
  \lambda & \nu \\
  \lambda^\dagger & \nu^\dagger
\end{array} \right]_{11}^{ij}\zeta_i(\kappa,\lambda^\dagger,\nu^\dagger)
\otimes\zeta_j(\kappa^\dagger,\nu,\lambda)\otimes
\zeta_l(\nu^\dagger,\mu,\lambda^\dagger)
\stackrel{\Id\otimes R \otimes \Id}{\longmapsto} \\
& & \sum_{\kappa, i,j,r} F_{0,\kappa}\left[ \begin{array}{cc}
  \lambda & \nu \\
  \lambda^\dagger & \nu^\dagger
\end{array} \right]_{11}^{ij}R_{jr}\zeta_i(\kappa,\lambda^\dagger,\nu^\dagger)
\otimes\zeta_r(\nu,\lambda,\kappa^\dagger)\otimes
\zeta_l(\nu^\dagger,\mu,\lambda^\dagger)
\stackrel{\Id\otimes F}{\longmapsto} \\
& & \sum_{\kappa, i,j,r} F_{0,\kappa}\left[ \begin{array}{cc}
  \lambda & \nu \\
  \lambda^\dagger & \nu^\dagger
\end{array} \right]_{11}^{ij} R_{jr}
F_{\nu,0}\left[ \begin{array}{cc}
  \lambda & \lambda^\dagger \\
  \kappa^\dagger & \mu
\end{array} \right]^{11}_{rl}
\zeta_i(\kappa,\lambda^\dagger,\nu^\dagger)
\otimes\zeta(\kappa^\dagger,\mu,0)\otimes
\zeta(0,\lambda^\dagger,\lambda) \stackrel{\Id\otimes R \otimes \Id}{\longmapsto} \\
& & \sum_{i,j,r} F_{0,\kappa}\left[ \begin{array}{cc}
  \lambda & \nu \\
  \lambda^\dagger & \nu^\dagger
\end{array} \right]_{11}^{ij} R_{jr}
F_{\nu,0}\left[ \begin{array}{cc}
  \lambda & \lambda^\dagger \\
  \kappa^\dagger & \mu
\end{array} \right]^{11}_{rl}
\zeta_i(\mu,\lambda^\dagger,\nu^\dagger)
\otimes\zeta(\mu^\dagger,\mu,0)\otimes
\zeta(0,\lambda^\dagger,\lambda) \stackrel{F \otimes \Id}{\longmapsto} \\
& & \sum_{i,j,r,s} F_{0,\kappa}\left[ \begin{array}{cc}
  \lambda & \nu \\
  \lambda^\dagger & \nu^\dagger
\end{array} \right]_{11}^{ij} R_{jr}
F_{\nu,0}\left[ \begin{array}{cc}
  \lambda & \lambda^\dagger \\
  \kappa^\dagger & \mu
\end{array} \right]^{11}_{rl}R_{is}
\zeta_s(\lambda^\dagger,\nu^\dagger,\mu)
\otimes\zeta(\lambda,0,\lambda^\dagger)\otimes
\zeta(0,\lambda^\dagger,\lambda) \stackrel{P^{(13)}(\Id\otimes R \otimes \Id)}{\longmapsto}\\
& & \sum_{i,j,r,s} F_{0,\kappa}\left[ \begin{array}{cc}
  \lambda & \nu \\
  \lambda^\dagger & \nu^\dagger
\end{array} \right]_{11}^{ij} R_{jr}
F_{\nu,0}\left[ \begin{array}{cc}
  \lambda & \lambda^\dagger \\
  \kappa^\dagger & \mu
\end{array} \right]^{11}_{rl}R_{is}
\zeta(0,\lambda^\dagger,\lambda)
\otimes\zeta(0,\lambda^\dagger,\lambda)\otimes
\zeta_s(\lambda^\dagger,\nu^\dagger,\mu) \stackrel{F \otimes \Id}{\longmapsto}\\
& & E_{\lambda^\dagger}\sum_{i,j,r,s} F_{0,\kappa}\left[ \begin{array}{cc}
  \lambda & \nu \\
  \lambda^\dagger & \nu^\dagger
\end{array} \right]_{11}^{ij} R_{jr}
F_{\nu,0}\left[ \begin{array}{cc}
  \lambda & \lambda^\dagger \\
  \kappa^\dagger & \mu
\end{array} \right]^{11}_{rl}R_{is}
\zeta(0,\lambda,\lambda^\dagger)
\otimes\zeta(0,\lambda,\lambda^\dagger)\otimes
\zeta_s(\lambda^\dagger,\nu^\dagger,\mu) \stackrel{\Id \otimes R \otimes \Id}{\longmapsto}\\
& & E_{\lambda^\dagger}\sum_{i,j,r,s} F_{0,\kappa}\left[ \begin{array}{cc}
  \lambda & \nu \\
  \lambda^\dagger & \nu^\dagger
\end{array} \right]_{11}^{ij} R_{jr}
F_{\nu,0}\left[ \begin{array}{cc}
  \lambda & \lambda^\dagger \\
  \kappa^\dagger & \mu
\end{array} \right]^{11}_{rl}R_{is}
\zeta(0,\lambda,\lambda^\dagger)
\otimes\zeta(\lambda,\lambda^\dagger,0)\otimes
\zeta_s(\lambda^\dagger,\nu^\dagger,\mu) \stackrel{\Id \otimes F}{\longmapsto}\\
& & E_{\lambda^\dagger}\sum_{i,j,r,s} F_{0,\kappa}\left[ \begin{array}{cc}
  \lambda & \nu \\
  \lambda^\dagger & \nu^\dagger
\end{array} \right]_{11}^{ij} R_{jr}
F_{\nu,0}\left[ \begin{array}{cc}
  \lambda & \lambda^\dagger \\
  \kappa^\dagger & \mu
\end{array} \right]^{11}_{rl}R^2_{is}
\zeta(0,\lambda,\lambda^\dagger)
\otimes\zeta(\nu, 0, \nu^\dagger)\otimes
\zeta_s(\nu^\dagger,\mu,\lambda^\dagger) \stackrel{\Id \otimes R \otimes \Id}{\longmapsto}\\
& & E_{\lambda^\dagger}\sum_{i,j,r,s} F_{0,\kappa}\left[ \begin{array}{cc}
  \lambda & \nu \\
  \lambda^\dagger & \nu^\dagger
\end{array} \right]_{11}^{ij} R_{jr}
F_{\nu,0}\left[ \begin{array}{cc}
  \lambda & \lambda^\dagger \\
  \kappa^\dagger & \mu
\end{array} \right]^{11}_{rl}R^2_{is}
\zeta(0,\lambda,\lambda^\dagger)
\otimes\zeta(0,\nu^\dagger,\nu)\otimes
\zeta_s(\nu^\dagger,\mu,\lambda^\dagger)
\end{eqnarray*}

This proves the stated formula.

\eproof

From formula (\ref{RelPent}) we conclude that $E_{\lambda} \neq 0$
and therefore also that
$S_{0,\lambda}\neq 0$ for all $\lambda \in \Lambda$.
Furthermore we also get the relation
\begin{equation}
 \sum_{i,j,r} F_{0,\mu}\left[ \begin{array}{cc}
  \lambda & \nu^\dagger \\
  \lambda^\dagger & \nu
\end{array} \right]_{11}^{lm} R_{lu}^2
F_{\nu^\dagger,0}\left[ \begin{array}{cc}
  \lambda & \lambda^\dagger \\
  \mu^\dagger & \mu
\end{array} \right]^{11}_{vu}R_{tv} = E_{\lambda^\dagger}^{-1}\delta_{tm}.\label{ABBA}
\end{equation}

By summing over $l$ in formula (\ref{RelPent}) we get
\begin{equation}
 \sum_{i,j,r} F_{0,\mu}\left[ \begin{array}{cc}
  \lambda & \nu \\
  \lambda^\dagger & \nu^\dagger
\end{array} \right]_{11}^{ij} R_{jr}F_{\nu,0}\left[ \begin{array}{cc}
  \lambda & \lambda^\dagger \\
  \mu^\dagger & \mu
\end{array} \right]^{11}_{rl}R^2_{is} = E_{\lambda^\dagger}^{-1}N_{\nu^\dagger,
\mu}^\lambda.\label{RelPentsum}
\end{equation}

As above, $\Theta$ is an oriented genus one surface and $\alpha$ and $\beta$ are simple closed
curves as indicated in figure \ref{figs-3}.
Let $\Thb'_\alpha$ and $\Thb'_\beta$ be marked surfaces, which
results from factoring $(\Si,L_\alpha)$ along $\alpha$, respectively
$(\Si,L_\beta)$ along $\beta$.
By factorization we get isomorphisms
\[\Phi'_\alpha : Z(\Theta,L_\alpha) \ra \bigoplus_{\lambda}
Z(\Thb'_\alpha,\lambda,\lambda^\dagger)\]
and
\[\Phi'_\beta : Z(\Theta,L_\beta) \ra \bigoplus_{\mu} Z(\Thb'_\beta,\mu,\mu^\dagger).\]

Pick diffeomorphisms
$$ f'_\alpha : \Xi \ra \Thb'_\alpha$$
and
$$ f'_\beta : \Xi \ra \Thb'_\beta$$
which maps the real axis onto $\beta$, respectively $\alpha$.

Then we get a basis $\xi_\lambda^\alpha$ for $Z(\Si,L_\alpha)$ and $\xi_\mu^\beta$
for $Z(\Si,L_\beta)$ by
\[\xi_\lambda^\alpha = (\Phi')^{-1}_\alpha Z(f'_\alpha)(\xi_\lambda)\]
and
\[\xi_\mu^\beta = (\Phi')^{-1}_\beta Z(f'_\beta)(\xi_\mu).\]

We have of course that
$$
\zeta_\lambda^\beta = \sum_\mu S_{\lambda,\mu}\zeta_\mu^\alpha.$$

\begin{proposition}\label{CurveOpformulav2}
We have the following formula
\[Z(\beta, \lambda)\xi_\mu^\alpha =
\sum_\nu N^\lambda_{\mu,\nu} \xi_\nu^\alpha.\]
\end{proposition}

\proof

We compute the action of $Z(\beta,\lambda)$ on $\zeta^\alpha_\mu$ by computing the
compositions indicated in
figure \ref{toruses2}:

\begin{center}
\begin{figure}
\includegraphics[scale=.5]{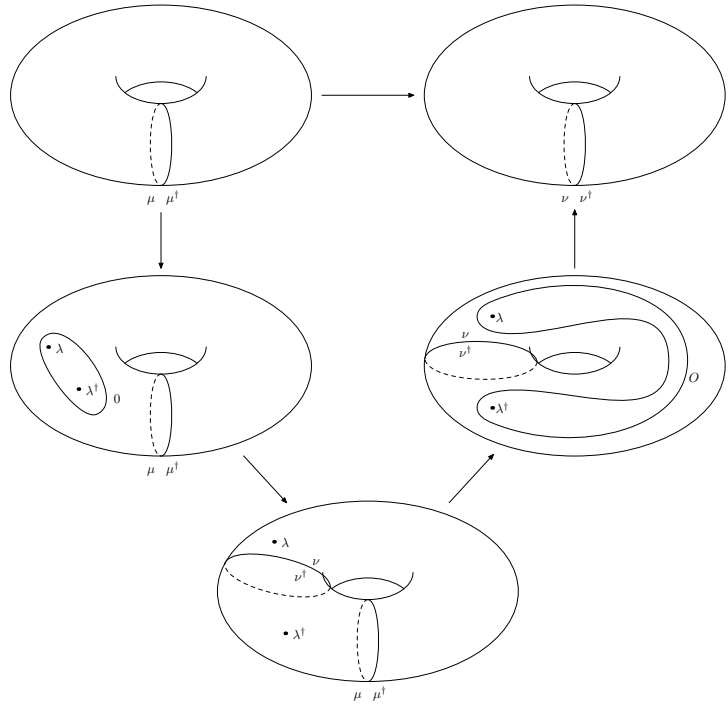}
\caption{The curve operator along $\beta$ on the genus one surface $\Theta$.}
\label{toruses2}
\end{figure}
\end{center}

\begin{eqnarray*}
\lefteqn{\zeta^\alpha_\mu \mapsto}\\
& & \zeta(0,\lambda,\lambda^\dagger) \otimes
\zeta(0,\mu,\mu^\dagger) \stackrel{F}{\longmapsto} \\
& & \sum_{\nu, i, j}F_{0,\nu}\left[ \begin{array}{cc}
  \lambda & \mu^\dagger \\
  \lambda^\dagger & \mu
\end{array} \right]_{11}^{ji}\zeta_j(\nu,\lambda^\dagger,\mu) \otimes
\zeta_i(\nu^\dagger,\mu^\dagger,\lambda) \stackrel{R^2\otimes R}{\longmapsto}  \\
& & \sum_{\nu, i, j,r,t}F_{0,\nu}\left[ \begin{array}{cc}
  \lambda & \mu^\dagger \\
  \lambda^\dagger & \mu
\end{array} \right]_{11}^{ji} R^2_{jr}R_{it}
\zeta_r(\mu,\nu,\lambda^\dagger) \otimes
\zeta_t(\mu^\dagger,\lambda,\nu^\dagger) \stackrel{F}{\longmapsto}  \\
& & \sum_{\nu, i, j,r,t}F_{0,\nu}\left[ \begin{array}{cc}
  \mu^\dagger & \lambda^\dagger \\
  \mu & \lambda
\end{array} \right]_{11}^{ji} R^2_{jr}R_{it}
F_{\mu^\dagger,0}\left[ \begin{array}{cc}
  \lambda & \lambda^\dagger \\
  \nu^\dagger & \nu
\end{array} \right]^{11}_{tr} \zeta(0,\nu^\dagger,\nu) \otimes
\zeta(0,\lambda^\dagger,\lambda) = \\
& & \sum_{\nu}E_{\lambda^\dagger}^{-1} N_{\mu, \nu}^\lambda \zeta(0,\nu^\dagger,\nu) \otimes
\zeta(0,\lambda^\dagger,\lambda) \mapsto
\sum_{\nu}E_{\lambda^\dagger}^{-1} N_{\mu, \nu}^\lambda \zeta^\alpha_\nu
\end{eqnarray*}
Hence
$$
Z(\beta,\lambda)(\zeta^\alpha_\mu) = \sum_\nu
N_{\mu,\nu}^{\lambda} \zeta_\nu^\alpha,$$
since
$$
\alpha_\lambda(Z(B) \zeta(\lambda,\lambda^\dagger,0)) =
E_{\lambda^\dagger}$$ by definition.
\eproof

\begin{proposition}\label{CurveOpformula}
We have the following formula
\[C_{\lambda,\mu} = S_{\mu,\lambda}/S_{\mu,0}.\]
\end{proposition}

\proof

From the formula in Proposition \ref{CurveOpformulav2} we deduce
that
\[S_{\mu,\rho} C_{\lambda,\mu} =  \sum_\nu N^\lambda_{\nu,\rho} S_{\mu,\nu}.\]
Letting $\rho = 0$, we get that
\[S_{\mu,0} C_{\lambda,\mu} = S_{\mu,\lambda},\]
which proves the stated formula, since we also conclude that
$S_{\mu,0} \neq 0$ from this.

\eproof

\begin{corollary}
We have the following formula for the coefficients $c_\kappa$:
\begin{equation}\label{cerne}
c_\kappa = d_{\kappa^\dagger} S_{\kappa^\dagger,0}.
\end{equation}
\end{corollary}

\proof One easily checks this formula by substitution. \eproof

\section{The formula for $S(\lambda)$.}
\label{sec8}

We begin by establishing a formula for the curve operator in terms
of the $F$, $R$, $B$ and the $d_\mu$'s.

\begin{lemma}\label{COF}
For all $\kappa,\lambda \in \Lambda$ and all
$j = 1, \ldots N_{\kappa, \nu}^{\nu}$ we have that
\begin{eqnarray*}
\lefteqn{Z(\beta, \kappa)(\zeta_i(\lambda,\mu,\mu^\dagger))}\\
 & = &
\sum_{\nu,j, r,k,m,s,w,p} d_\nu^{-1}R_{ip} B_{pr}F_{\nu^\dagger,\mu}\left[ \begin{array}{cc}
  \kappa & \mu^\dagger \\
  \nu^\dagger & \lambda
\end{array} \right]_{kr}^{sm} R_{mk} R^2_{sw} B_{wj}
\zeta_j(\kappa,\nu,\nu^\dagger).
\end{eqnarray*}
\end{lemma}

\begin{center}
\begin{figure}
\includegraphics[scale=.5]{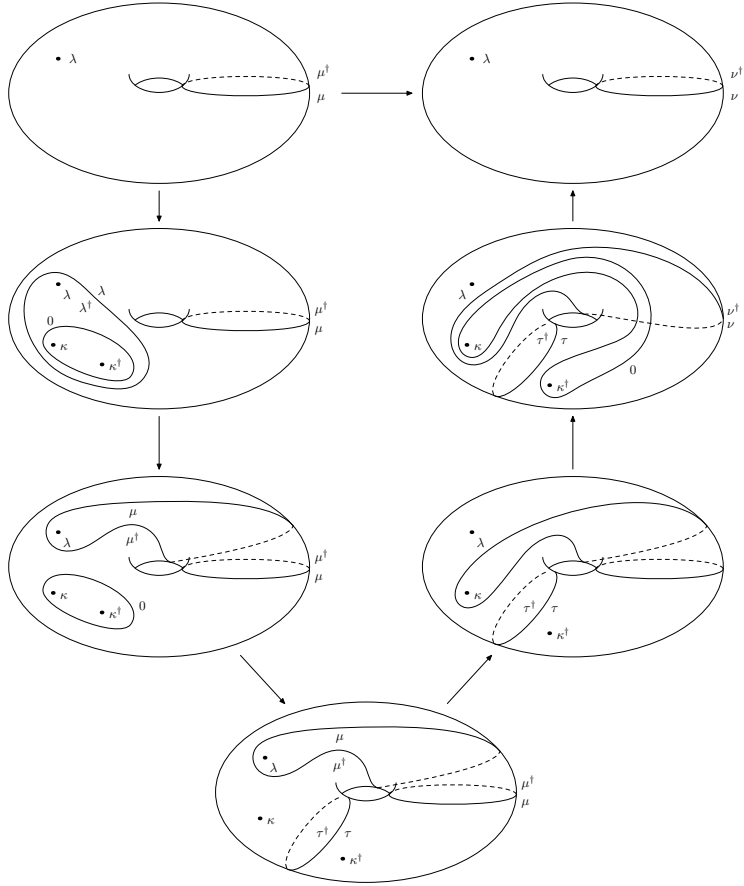}
\caption{The curve operator along $\beta$ via changes $\M$ to $\cF$.}
\label{COviaFmoves}
\end{figure}
\end{center}

\proof

The curve operator $Z(\beta,\lambda)$ acts as follows:

\begin{eqnarray*}
\lefteqn{\zeta_i(\lambda,\mu,\mu^\dagger) \mapsto}\\
& &
\zeta(\lambda^\dagger,\lambda,0) \otimes \zeta_i(\lambda,\mu,\mu^\dagger)
\otimes
\zeta(0,\kappa,\kappa^\dagger)
\stackrel{(\Id\otimes R \otimes \Id)P^{231}(F\otimes \Id)}{\longmapsto}\\
& &
\sum_p R_{ip}\zeta(0,\kappa,\kappa^\dagger) \otimes
\zeta(0,\mu,\mu^\dagger) \otimes
\zeta_p(\mu,\mu^\dagger,\lambda)
 \stackrel{F\otimes \Id}{\longmapsto}\\
& & \sum_{\tau, p,l, m}R_{ip} F_{0,\tau}\left[ \begin{array}{cc}
  \kappa & \mu^\dagger \\
  \kappa^\dagger & \mu
\end{array} \right]_{11}^{lm}\zeta_l(\tau,\kappa^\dagger,\mu)\otimes
\zeta_m(\tau^\dagger,\mu^\dagger,\kappa)\otimes
\zeta_p(\mu,\mu^\dagger,\lambda)
\stackrel{\Id\otimes R \otimes B_{23}}{\longmapsto} \\
& & \sum_{\tau, p,k,l, m,r}R_{ip} (B_{23})_{pr} R_{mk}
F_{0,\tau}\left[ \begin{array}{cc}
  \kappa & \mu^\dagger \\
  \kappa^\dagger & \mu
\end{array} \right]_{11}^{lm}\zeta_l(\tau,\kappa^\dagger,\mu)\otimes
\zeta_k(\mu^\dagger,\kappa,\tau^\dagger)\otimes
\zeta_r(\mu,\lambda,\mu^\dagger)  \stackrel{\Id\otimes F}{\longmapsto} \\
& & \sum_{\tau,\nu, l,m,p,k,r,s,t}R_{ip} (B_{23})_{pr} R_{mk}
F_{0,\tau}\left[ \begin{array}{cc}
  \kappa & \mu^\dagger \\
  \kappa^\dagger & \mu
\end{array} \right]_{11}^{lm}
F_{\mu\dagger,\nu}\left[ \begin{array}{cc}
  \kappa &  \mu^\dagger\\
  \tau^\dagger & \lambda
\end{array} \right]^{st}_{kr}\\
& & \phantom{hejhej}\zeta_l(\tau,\kappa^\dagger,\mu)\otimes
\zeta_s(\nu,\tau^\dagger,\lambda)\otimes
\zeta_t(\nu^\dagger,\mu^\dagger,\kappa)
\stackrel{P^{231}(R^2\otimes \Id \otimes R)}{\longmapsto} \\
& & \sum_{\tau,\nu, l,m,p,k,r,s,t,u,v}R_{ip} (B_{23})_{pr} R_{mk}
F_{0,\tau}\left[ \begin{array}{cc}
  \kappa & \mu^\dagger \\
  \kappa^\dagger & \mu
\end{array} \right]_{11}^{lm}
F_{\mu\dagger,\nu}\left[ \begin{array}{cc}
  \kappa &  \mu^\dagger\\
  \tau^\dagger & \lambda
\end{array} \right]^{st}_{kr} R^2_{lu} R_{tv}\\
& & \phantom{hejhej}
\zeta_s(\nu,\tau^\dagger,\lambda)\otimes
\zeta_v(\mu^\dagger,\kappa,\nu^\dagger)\otimes
\zeta_u(\mu,\tau,\kappa^\dagger)
\stackrel{\Id \otimes F}{\longmapsto} \\
& & \sum_{\nu, l,m,p,k,r,s,t,u,v}R_{ip} (B_{23})_{pr} R_{mk}
F_{0,\nu}\left[ \begin{array}{cc}
  \kappa & \mu^\dagger \\
  \kappa^\dagger & \mu
\end{array} \right]_{11}^{lm}
F_{\mu\dagger,\nu}\left[ \begin{array}{cc}
  \kappa &  \mu^\dagger\\
  \nu^\dagger & \lambda
\end{array} \right]^{st}_{kr} R^2_{lu} R_{tv}\\
& & \phantom{HejHej} F_{\mu^\dagger,0}\left[ \begin{array}{cc}
  \kappa & \kappa^\dagger \\
  \nu^\dagger & \nu
\end{array} \right]_{vu}^{11}
\zeta_s(\nu,\nu^\dagger,\lambda)\otimes
\zeta(0,\nu^\dagger,\nu)\otimes
\zeta(0,\kappa^\dagger,\kappa)
\end{eqnarray*}
We now apply formula (\ref{ABBA}) to this expression and we get that
\begin{eqnarray*}
\lefteqn{\zeta_i(\lambda,\mu,\mu^\dagger) \mapsto}\\
& & \sum_{\nu, l,m,p,k,r,s}R_{ip}  (B_{23})_{pr}
F_{\mu\dagger,\nu}\left[ \begin{array}{cc}
  \kappa &  \mu^\dagger\\
  \nu^\dagger & \lambda
\end{array} \right]^{sm}_{kr}  R_{mk} \\
& & \phantom{HejHej}
\zeta_s(\nu,\nu^\dagger,\lambda)\otimes
\zeta(0,\nu^\dagger,\nu)\otimes
\zeta(0,\kappa^\dagger,\kappa) \mapsto\\
& & \sum_{\nu, l,m,p,k,r,s}R_{ip}  (B_{23})_{pr}
F_{\mu\dagger,\nu}\left[ \begin{array}{cc}
  \kappa &  \mu^\dagger\\
  \nu^\dagger & \lambda
\end{array} \right]^{sm}_{kr}  R_{mk} \\
& & \phantom{HejHej}
(B_{12})_{sw} R^2_{wj} d_\mu^{-1}\zeta_j(\lambda,\nu,\nu^\dagger)
\end{eqnarray*}
from which the formula follows.
\eproof

{\noindent {\bf Proof of Theorem \ref{Main'}.}}
By the calculations in section \ref{From1to0} and Proposition \ref{CODTp}
we get that
$$
S(\lambda)_{\mu,i}^{\nu,j} = d_\nu \sum_\kappa c_\kappa
Z(\beta,\kappa)_{\mu,i}^{\nu,j} d_\mu.$$

Combining this with formula (\ref{cerne}) and the formula for
$Z(\beta,\kappa)$ from Proposition \ref{COF}, we obtain
the formula stated in Theorem \ref{Main'}.

\eproof

\end{document}